\documentclass[sn-apa]{sn-jnl}

\usepackage{graphicx}%
\usepackage{multirow}%
\usepackage{amsmath,amssymb,amsfonts}%
\usepackage{amsthm}%
\usepackage{mathrsfs}%
\usepackage[title]{appendix}%
\usepackage{xcolor}%
\usepackage{textcomp}%
\usepackage{manyfoot}%
\usepackage{booktabs}%
\usepackage{algorithm}%
\usepackage{algorithmicx}%
\usepackage{algpseudocode}%
\usepackage{listings}%
\usepackage{longtable}
\usepackage{subcaption}
 \usepackage{geometry}
\usepackage{orcidlink}
\theoremstyle{thmstyleone}%
\theoremstyle{thmstyletwo}%

\theoremstyle{thmstylethree}%

\begin{document}

\title[Two-stage heuristic algorithm for a new variant of the multi-compartment vehicle routing problem with stochastic demands]{Two-stage heuristic algorithm for a new variant of the multi-compartment vehicle routing problem with stochastic demands}


\author[1]{\fnm{Juan Carlos} \sur{Gon\c calves-Dosantos}\orcidlink{0000-0003-1578-8411}}\email{jgoncalves@umh.es}
\equalcont{These authors contributed equally to this work.}

\author*[2,3]{\fnm{Laura} \sur{Davila-Pena}\orcidlink{0000-0003-2175-2546}}\email{l.davila-pena@kent.ac.uk; lauradavila.pena@usc.es}
\equalcont{These authors contributed equally to this work.}

\author[4]{\fnm{Balbina} \sur{Casas-Méndez}\orcidlink{0000-0002-2826-218X}}\email{balbina.casas.mendez@usc.es}
\equalcont{These authors contributed equally to this work.}

\affil[1]{\orgname{University Miguel Hern{\'a}ndez de Elche}, \city{Elche}, \postcode{03202},  \country{Spain}}

\affil*[2]{\orgdiv{Centre for Sustainable Logistics and Analytics (CeLSA), Department of Analytics, Operations and Systems}, \orgname{Kent Business School, University of Kent}, \orgaddress{\street{Canterbury Campus}, \city{Canterbury}, \postcode{CT2 7PE}, \country{UK}}}

\affil[3]{\orgdiv{MODESTYA Research Group, Department of Statistics, Mathematical Analysis and Optimization, Faculty of Mathematics}, \orgname{University of Santiago de Compostela}, \orgaddress{\street{Campus Vida}, \city{Santiago de Compostela}, \postcode{15782}, \country{Spain}}}

\affil[4]{\orgdiv{CITMAga, MODESTYA Research Group, Department of Statistics, Mathematical Analysis and Optimization, Faculty of Mathematics}, \orgname{University of Santiago de Compostela}, \orgaddress{\street{Campus Vida}, \city{Santiago de Compostela}, \postcode{15782}, \country{Spain}}}


\abstract{This paper presents a model for a vehicle routing problem in which customer demands are stochastic and vehicles are divided into compartments. The problem is motivated by the needs of certain agricultural cooperatives that produce various types of livestock food. The vehicles and their compartments have different capacities, and each compartment can only contain one type of feed. Additionally, certain farms can only be accessed by specific vehicles, and there may be urgency constraints. To solve the problem, a two-step heuristic algorithm is proposed. First, a constructive heuristic is applied, followed by an improvement phase based on iterated tabu search. The designed algorithm is tested on several instances, including an analysis of real-world datasets where the results are compared with those provided by the model. Furthermore, multiple benchmark instances are created for this problem and an extensive simulation study is conducted. Results are presented for different model parameters, and it is shown that, despite the problems' complexity, the algorithm performs efficiently. Finally, the proposed heuristic is compared to existing solution algorithms for similar problems using benchmark instances from the literature, achieving competitive results.}




\maketitle

\section{Introduction}\label{sec:intro}

Vehicle routing problems (VRPs) aim to design a set of minimum-cost routes for a fleet of vehicles to meet the demands of a group of customers \citep{Dan59}. More complex variants of the classic VRP emerge when additional constraints are introduced to better reflect real-world settings (see, for example, \citealp{Fug2006}; \citealp{Rev2017};  \citealp{Grab2021}). One such variant is the multi-compartment vehicle routing problem (MC-VRP; \citealp{Bro81} and \citealp{Bro87}). While MC-VRPs share the same objective as the classic VRP, they involve transporting different, often incompatible, products that must be carried in separate compartments within the vehicles. Seminal formulations for these problems can be found in \cite{elf08} and \cite{Der13}, with \cite{Coe2015} also offering a classification scheme for MC-VRPs. The use of multiple compartments was first motivated by the petrol station replenishment problem \citep{Corn2008}. Since then, significant research has been conducted on the petroleum distribution problem in the context of MC-VRPs, with recent studies including \cite{Sun2021} and \cite{Haj2024}. Other common applications for multi-compartmentalized vehicles include waste collection \citep{Hess2024} and food delivery \citep{Dav23}. For extensive reviews of different applications, we refer the reader to \cite{Der2011}, \cite{Ost2021}, and \cite{Gu2024}. 


Most of these works, however, overlook an important aspect: they fail to account for the randomness of real-life parameters. Even with a large number of constraints, these models do not adequately address the inherent uncertainty in real-world scenarios. When any of the model parameters are random, the problem becomes a stochastic vehicle routing problem (SVRP). Unlike deterministic models, SVRPs have a less cohesive and more scattered literature (\citealp{Oyo17, Oyo18}). One of the most extensively studied variants in this context is the vehicle routing problem with stochastic demands (VRPSD; \citealp{Ber92}), where customer demands are treated as random variables. As in the deterministic case, the goal is to determine a set of routes that minimizes the total distance, which includes both the distance traveled on the fixed routes and the expected extra distance--referred to as the ``expected distance''--resulting from fluctuations in demand. These expected distances occur when the actual demand exceeds the planned capacity, forcing the vehicle to return to the depot mid-route before continuing. Since the demand is known gradually as the vehicle visits customers, route adjustments may be necessary during operation. As noted by \cite{Ber92}, redesigning the routes once the demand becomes known may require a big effort and cost. Additionally, the route operator may prioritize service consistency and customization, such as assigning the same driver and vehicle to a specific customer. However, the exact demand for such a customer remains unknown until the visit takes place. These practical considerations, along with a real-world case study, motivate our stochastic approach.  

This work integrates the two aspects of VRPs mentioned above: multi-compartment vehicles and stochastic demands, resulting in a multi-compartment vehicle routing problem with stochastic demands (MC-VRPSD). There is limited research in the literature addressing this combined problem, which was first introduced in \cite{Men10,Men11}. Building on these early contributions, we incorporate several constraints tailored to the specific needs of a Spanish agricultural cooperative that produces and distributes various types of livestock food. These constraints include prohibiting the same compartment from carrying different products or products for different customers (due to the vehicles' inability to measure partial discharges within a compartment), handling urgent orders (a type of so-called time window constraint), and accounting for potential accessibility restrictions for trucks at certain customer locations. In addition, the cooperative is interested not only in minimizing the transportation costs, but also in maximizing the transported load. Collectively, all these elements contribute to the novelty of the model we are developing, to the best of our knowledge.

Determining the optimal routes for the MC-VRPSD has an NP-hard complexity. While exact algorithms can be useful for testing the model, they are generally practical only for instances with a small number of customers. For larger cases, metaheuristic algorithms, often combined with heuristics that generate a good initial solution, are typically considered for these VRP variants. For instance, \cite{Muy10} present a three-stage algorithm for the MC-VRP, where an initial feasible solution is obtained using the Clarke and Wright savings algorithm \citep{Cla64}, which is then refined by local search moves, all integrated within the guided local search metaheuristic. \cite{Yah2020} introduce two metaheuristics to solve MC-VRPs: an adaptive variable neighborhood search and a genetic algorithm. Tabu-search based algorithms are particularly popular for MC-VRPs, with \cite{Sil17} and \cite{Dav23} developing iterated tabu search metaheuristics, while \cite{elf08} explore both a tabu search and a memetic algorithm. A memetic algorithm is also introduced in \cite{Men10} to solve the MC-VRPSD. In \cite{Men11}, the authors propose a set of construction heuristics, including stochastic versions of the nearest neighbor, best insertion, and \cite{Dror86} savings algorithm, all adapted to the multi-compartment case.

Considering all the aspects discussed, we have tested our proposed model on small-size instances and developed a two-stage heuristic approach. In the first stage, we designed a constructive algorithm that extends the well-known Clarke and Wright savings algorithm \citep{Cla64}. When determining the savings associated with inserting a new customer into a route, this initial solution incorporates constraints that prevent the mixing of products within the same truck compartment and considers urgent orders. Additionally, recognizing the stochastic nature of the problem and the potential for high demand in relation to compartment capacity, we allow a customer's demand to be divided and served across different compartments of the same truck or different trucks. Our proposed final solution is an improvement over this initial approach, achieved through a metaheuristic method based on an iterated tabu search that employs strategies for the destruction and recreation of routes. To validate the model, we use the \texttt{AMPL} (\url{https://ampl.com/}) algebraic modeling language and the \texttt{Gurobi} (\url{https://www.gurobi.com/}) solver, while the algorithms are implemented in the \texttt{R} programming language (\url{https://cran.r-project.org/}). We test the algorithms using real-life data and also create a set of test cases for a comprehensive computational study. In both scenarios, we present results based on varying parameter values of the model. Since the MC-VRPSD model for which the two-stage algorithm was designed is novel, direct comparisons with existing methods in the literature are not possible. However, this limitation is mitigated by applying the algorithm to alternative models, allowing for comparisons with exact solutions or other algorithms. Our results indicate that our algorithm is competitive in achieving near-optimal solutions. In particular, we analyzed its performance on a small real-world problem, for which the exact solution is known, while considering various fleet sizes and scenarios with and without stochasticity. Additionally, we examined the algorithm's performance in the MC-VRP context using instances from the literature and the results provided in \cite{elf08} and \cite{Men10}.

The remainder of this paper is organized as follows. Section~\ref{sec:mcvrpsd} contextualizes the problem and introduces the mathematical formulation. Section~\ref{sec:heuristic} describes the proposed two-phase heuristic algorithm. In Section~\ref{sec:results}, the results of the computational study are discussed, including the exact solving of small-sized problems, the application of the proposed heuristic to existing and newly created benchmark problems, and the comparison with previous methods in the literature. Finally, concluding remarks are presented in Section~\ref{sec:conclusions}.

\section{Multi-compartment vehicle routing problem with stochastic demands (MC-VRPSD)}\label{sec:mcvrpsd}

\subsection{Contextual description of the problem}\label{sec:description}

Certain agricultural cooperatives produce various types of animal feed, which are distributed to numerous farmers across a broad area. Customers often place orders for different types of feed, and the cooperative has specific delivery times based on the urgency of these orders. Urgent orders generally arise when a customer has a limited amount of feed, requiring a prompt refill of their tank to ensure that their animals have enough food. If an order is classified as urgent, it must be delivered within one day, resulting in higher delivery costs. This practice can lead to system saturation, which is not beneficial for efficient route design.  

In this paper, we focus on a cooperative located in Galicia, a region in the northwest of Spain, which produces four types of animal feed. This cooperative serves over 1500 customers, each of whom typically places up to two orders per month. The cooperative operates a fleet of different trucks, each equipped with multiple compartments of varying capacities. Each compartment can only transport one type of feed and can accommodate only feed for one customer. The trucks are also subject to restrictions on the distance traveled and the load carried each day. Drivers are compensated based on the distance traveled, which can also be expressed in terms of the time taken for delivery, as well as the amount of cargo transported. Additionally, access to individual farms may be limited to certain vehicles. 

The primary objective of this work is to provide the cooperative with an anticipation tool--a method for designing routes for future deliveries without waiting for customer orders. This approach aims to improve the distribution efficiency, reduce associated costs and effort, and, most importantly, prevent system saturation. By knowing the number of days since the last delivery and the estimated daily feed consumption of each farm, we can assess the urgency of the next visit. Furthermore, historical order data enables us to estimate the probabilistic distribution of the customer demand. Moreover, advances in technology also allow for monitoring the inventory levels of cooperative customers. With this information, the cooperative can select which customers to serve each day according to the available fleet. However, there is a need for a tool that effectively plans vehicle routes to maximize transported load and minimize transportation costs. To address this problem comprehensively, our model will account for a heterogeneous fleet of vehicles, considering both their capacities and the number and capacity of their compartments.

\subsection{Mathematical formulation of the MC-VRPSD}\label{sec:formulation}

In this section, we present the proposed mathematical formulation to model our MC-VRPSD. 

\subsubsection{Indices and sets}\label{sec:index_sets}
The main elements of our formulation are the following:

\begin{longtable*}{llp{0.9\textwidth}}

$n,N$ &= &index and set of customers,\\[0.15cm]
  
$\bar{N}:=N\cup \{0\}$ &=& set of nodes, including customers and the cooperative headquarters (depot),\\[0.15cm]
  
$k,K$&= &index and set of trucks,\\[0.15cm]

$h,H_k$&=& index and set of compartments for truck $k\in K$,\\[0.15cm]

$f,F$&=& index and set of livestock feeds,\\[0.15cm]

$r,R_k$&= &index and set of possible routes for truck $k\in K$.
\end{longtable*}

\subsubsection{Parameters}\label{sec:parameters}

Here we present the parameters that will be included in the model, encompassing both deterministic and stochastic elements.\\
 
\begin{longtable*}{llp{0.7\textwidth}}

$\displaystyle C^{k}:=\{c_{h}^{k}\}_{h\in H_k}$&=&  vector of capacities of truck $k\in K$, where $c_{h}^{k}$ is the capacity of compartment $h\in H_k$,\\[0.15cm]

$\displaystyle \textbf{O}:=\left[\textbf{O}_{n,f}\right]_{n\in N, f\in F}$&=&  matrix of stochastic orders, where $\textbf{O}_{n,f}$ is the  random variable describing the non-negative order of customer $n\in N$ for feed $f\in F$,\\[0.15cm]

$\displaystyle D:=\left[d_{n_{1},n_{2}}\right]_{n_1,n_2 \in \bar{N}}$&=&  matrix of distances, in minutes, where $d_{n_{1},n_{2}}$ is the  time a truck needs to go from node $n_1\in \bar{N}$ to node $n_2\in \bar{N}$,\\[0.15cm]

$\displaystyle I:=\left[i_{k,n}\right]_{k\in K,n\in N}$&=&  matrix that describes which customers can be visited by truck $k\in K$, such that $i_{k,n}$ equals $1$ if truck $k$ can visit customer $n\in N$, and $0$ otherwise, \\[0.15cm]

$l_{k}$&=&  maximum allowable load of truck $k\in K$,\\[0.15cm]

$\displaystyle \textbf{U}:= \left[\textbf{U}_{n,f}\right]_{n\in N, f\in F}$&=&  matrix of stochastic urgency, where $\textbf{U}_{n, f}$ is the random variable describing the maximum number of days that customer $n\in N$ can wait to receive feed $f\in F$.
\end{longtable*}

\subsubsection{Decision variables}\label{sec:variables}

The following variables define the routes traveled by the trucks and how their compartments are loaded.\\

\begin{longtable*}{llp{0.75\linewidth}}
    $x_{n_{1},n_{2}}^{k,r} $  & $=$ &  $\left\lbrace
\begin{array}{ll}
1 & {\rm if}\ {\rm truck}\ k\in K, \ {\rm on}\ {\rm route}\ r\in R_k,\ {\rm travels}\ {\rm from}\ n_1\in \bar{N}\ {\rm to}\ 
  n_2\in \bar{N},\\
0 & {\rm otherwise}.
\end{array}
\right.$\\[0.15cm]

$y_{n,f,h}^{k,r} $ & $=$ &  Proportion of compartment $h\in H_{k}$ used to transport feed $f\in F$ by truck $k\in K$ on route $r\in R_{k}$ to serve customer $n\in N$.\\[0.15cm]

$z_{n,f}^{k,r}$ & $=$ &  Probability that an urgent demand for feed $f\in F$ from  customer $n\in N$ exceeds the quantity delivered by truck $k\in K$ on route $r\in R_k$.\\[0.15cm]

$v_{n,f,h}^{k,r} $ & $=$ & $\left\lbrace
\begin{array}{ll}
1 & \textrm{if compartment}\ h\in H_{k}\ \textrm{is used to carry feed} \ f\in F\ \textrm{by truck}\ k\in K, \textrm{on}\ \\& \textrm{route}\ r\in R_{k} \  \textrm{to serve customer}\ n\in N,\\
0 & {\rm otherwise}.
\end{array}
\right.$\\[0.15cm]

$u_{n_1,n_2,f_1,f_2,h}^{k,r}$ & $=$ & $\max\{y_{n_1,f_1,h}^{k,r},y_{n_2,f_2,h}^{k,r}\}$, where $h\in H_{k}$, $f_1, f_2\in F$, $k\in K$, $r\in R_{k}$, and $n_1, n_2\in N$.
\end{longtable*}

\subsubsection{Constraints}\label{sec:constraints}

The constraints of the mathematical model are classified into five groups, which we will explain below.

\begin{enumerate}
\item[1)] {Constraints describing the truck routes.
\begin{equation}
x^{k,r}_{n,n}=0 \hspace*{6.0cm}  \forall k\in K, \forall r\in R_{k}, \forall n\in \bar{N}.
\label{1.1}
\end{equation}
\begin{equation}
\sum_{n\in N}x_{0,n}^{k,r}\leq1 \hspace*{6.7cm} \forall k\in K, \forall r\in R_{k}.
\label{1.2}
\end{equation}
\begin{equation}
\displaystyle x^{k,r}_{n_{1},n_{2}}\leq\sum_{n\in N}x_{0,n}^{k,r} \hspace*{4.0cm} \forall k\in K, \forall r\in R_{k}, \forall n_{1},n_{2}\in N.
\label{1.3}
\end{equation}
\begin{equation}
\displaystyle\sum_{n\in \bar{N}}x_{n,n_{1}}^{k,r}\leq\sum_{n\in \bar{N}}x_{n_{1},n}^{k,r} \hspace*{3.7cm} \forall k\in K, \forall r\in R_{k}, \forall n_{1}\in \bar{N}.
\label{1.4}
\end{equation}
\begin{equation}
\displaystyle\sum_{n_1,n_2\in S}x_{n_{1},n_{2}}^{k,r}\leq|S|-1 
\hspace*{2.0cm} \forall k\in K, \forall r\in R_{k}, \forall S\subseteq N\ {\rm and}\ |S|>1.
\label{1.5}
\end{equation}
\vspace*{0.05cm}

Constraint~\eqref{1.1} ensures that no truck can travel from a customer back to the same customer. Constraint~\eqref{1.2} limits each truck route to departing from the cooperative headquarters at most once. Constraint~\eqref{1.3} requires all planned routes to start at the cooperative headquarters. Constraint~\eqref{1.4} is a flow conservation constraint, meaning that if a truck arrives at a customer, it must also leave from that customer. Lastly, Constraint~\eqref{1.5} prevents the formation of sub-tours on any truck route.\\
    }

\item[2)] {Constraints linking the truck routes with their loading.
\begin{equation}
\displaystyle \frac{1}{|H_{k}|} \sum_{f\in F}\sum_{h\in H_{k}}y_{n_{1},f,h}^{k,r}\leq\sum_{n\in \bar{N}}x_{n,n_{1}}^{k,r} \hspace*{2.1cm} \forall k\in K, \forall r\in R_{k}, \forall n_{1}\in N.
\label{2.1}
\end{equation}
\begin{equation}
\displaystyle \sum_{n\in \bar{N}}x_{n,n_{1}}^{k,r}\leq\sum_{f\in F}\sum_{h\in H_{k}}M'  y_{n_{1},f,h}^{k,r} \hspace*{2.4cm} \forall k\in K, \forall r\in R_{k}, \forall n_{1}\in N.
\label{2.2}
\end{equation}
\vspace*{0.05cm}

Constraint~\eqref{2.1} ensures that if a truck is carrying feed for a customer, it must visit that customer. Constraint~\eqref{2.2} guarantees that if a truck visits a customer on a route, it is carrying feed for them. Here, $M'$ represents a sufficiently large constant. In general, it will be enough to take $M'=c_{h}^{k}$, for some $h\in H_k$.\\

    }

    \item[3)] {Constraints modeling the random urgency of customer orders.
\begin{equation}
\hspace{-0.5cm}  \mathbb{P}\left[\sum_{k\in K}\sum_{r\in R_{k}}\sum_{h\in H_{k}}c_{h}^{k}y_{n,f,h}^{k,r}\geq \textbf{O}_{n,f}\right]\leq \mathbb{P}\left[ \textbf{U}_{n,f}=0\right] \hspace*{1.1cm} \forall f\in F, \forall n\in N.
\label{3.1}
\end{equation}
\begin{equation}
\begin{split}
\hspace{1.4cm} &\mathbb{P}\left[\sum_{k\in K}\sum_{r\in R_{k}}\sum_{h\in H_{k}}c_{h}^{k}y_{n,f,h}^{k,r}\geq \textbf{O}_{n,f}\right]> \epsilon\ (\approx 0) \hspace*{1.8cm} \forall n\in \hat{N},  \forall f\in \hat{F}_{n}, \label{3.2}\\&
\textrm{where}\ \hat{N}\subseteq N\ \textrm{and}\ \forall n\in\hat{N}, \hat{F}_n=\{f\in F\,:\, \mathbb{P}\left[ \textbf{U}_{n,f}=0\right]\geq \beta\}\ (\textrm{such that}\ \beta=0.90).
\end{split}
\end{equation}
\vspace*{0.05cm}

Constraint~\eqref{3.1} indicates that as the probability of a demand being very urgent decreases, the probability of fully satisfying that demand also decreases. Constraint~\eqref{3.2} states that there is a positive probability of meeting the demands of customers with a high probability of urgent needs. $\mathbb{P}$ represents the probability function. Additionally, $\hat{N}$ refers to the set of customers for whom the demand for at least one type of feed is urgent with a high probability. For each customer $n\in \hat{N}$, the corresponding set of feeds (whose demand is urgent with a high probability) is denoted by $\hat{F}$.\\
    }

\item[4)]{Constraints modeling accessibility, truck capacities, and limitations on time (minutes)  
during the working day (capacity constraints).

\begin{equation}
\displaystyle \sum_{n\in \bar{N}} x_{n,n_{1}}^{k,r}\leq i_{k,n_{1}}  \hspace*{4.4cm} \forall k\in K, \forall r\in R_{k}, \forall n_{1}\in N.
\label{4.1}
\end{equation}
\begin{equation}
\displaystyle \sum_{f\in F} y_{n_{1},f,h}^{k,r}\leq i_{k,n_{1}}  \hspace*{2.8cm} \forall k\in K, \forall r\in R_{k}, \forall h\in H_{k} , \forall n_{1}\in N.
\label{4.2}
\end{equation}
\begin{equation}
\displaystyle \sum_{r\in R_{k}}\sum_{n_{1}\in \bar{N}}\sum_{n_{2}\in \bar{N}} d_{n_{1},n_{2}}x_{n_{1},n_{2}}^{k,r}  +
\sum_{r\in R_{k}}\sum_{n\in \hat{N}}\sum_{f\in \hat{F}_n}2 d_{0,n}z_{n,f}^{k,r}\leq 540 \hspace*{0.9cm} \forall k\in K.
\label{4.3}
\end{equation}
\begin{equation}
\begin{split}
& z_{n,f}^{k,r}\leq p_{n,f}^{k,r} \hspace*{4.2cm} \forall k\in K,  \forall r\in R_{k}, \forall n\in\hat{N}, \forall f\in \hat{F}_n,\\&
\textrm{where}\ p_{n,f}^{k,r}=\mathbb{P}\left[\displaystyle \sum_{h\in H_k} c_{h}^{k}y_{n,f,h}^{k,r}< \textbf{O}_{n,f}\right]. 
\end{split}
\label{4.3.2}
\end{equation}
\begin{equation}
\displaystyle z_{n,f}^{k,r}\leq \sum_{n_1\in \bar{N}}x_{n_{1},n}^{k,r} \hspace*{3.10cm} \forall k\in K,  \forall r\in R_{k}, \forall n\in\hat{N}, \forall f\in \hat{F}_n.
\label{4.3.3}
\end{equation}
\begin{equation}
\displaystyle z_{n,f}^{k,r}\geq (\sum_{n_1\in \bar{N}}x_{n_{1},n}^{k,r}-1)+  p_{n,f}^{k,r} \hspace*{1.40cm} \forall k\in K,  \forall r\in R_{k},  \forall n\in\hat{N}, \forall f\in \hat{F}_n.
\label{4.3.4}
\end{equation}
\begin{equation}
\displaystyle \sum_{n\in N}\sum_{f\in F}\sum_{h\in H_{k}}c_{h}^{k}y_{n,f,h}^{k,r}\leq l_{k} \hspace*{4.5cm} \forall k\in K, \forall r\in R_{k}.
\label{4.5}
\end{equation}
\vspace*{0.05cm}

Constraints~\eqref{4.1} and \eqref{4.2} ensure that impossible visits are not allowed. Constraint~\eqref{4.3} establishes that each truck route has a maximum duration of  $540$ minutes. Note that we do not consider unloading time. To accurately formulate this restriction, we need to sum the expected time for what is referred to as a failed route (the second term on the left side of the constraint). This occurs when a truck visits a customer with an urgent order but lacks sufficient cargo to meet the full demand, necessitating a return to the depot. 
Constraints~\eqref{4.3.2}--\eqref{4.3.4} ensure that $z_{n,f}^{k,r}$ takes indeed its proper value, i.e, is equal to the probability that an urgent demand for feed $f$ from customer $n$ exceeds the quantity delivered by truck $k$ on route $r$, provided that the customer is part of that route.
Lastly, Constraint~\eqref{4.5} indicates that a truck cannot transport more than its legally authorized cargo.\\
}

\item[5)] {Constraints modeling technical limitations on the truck loading procedure.
\begin{equation}
\begin{array}{l}
\displaystyle y_{n_{1},f_{1},h}^{k,r}+y_{n_{2},f_{2},h}^{k,r}\leq u^{k,r}_{n_1,n_2,f_1,f_2,h} \hspace*{2.5cm}\forall k\in K, \forall r\in R_{k}, \forall h\in H_{k},\\
\hspace*{7.5cm} \forall n_{1},n_{2}\in N, \forall f_{1},f_{2}\in F.
\end{array}\label{17}
\end{equation}
\begin{equation}
\begin{array}{l}
\displaystyle u^{k,r}_{n_1,n_2,f_1,f_2,h}\leq y_{n_{1},f_{1},h}^{k,r}+  10\cdot(1-  v^{k,r}_{n_1,f_1,h})  \hspace*{1cm}\forall k\in K, \forall r\in R_{k}, \forall h\in H_{k},\\
\hspace*{7.5cm} \forall n_{1}, n_2\in N, \forall f_{1},f_{2}\in F.
\end{array}\label{18}
\end{equation}
\begin{equation}
\begin{array}{l}
\displaystyle u^{k,r}_{n_1,n_2,f_1,f_2,h}\leq y_{n_{2},f_{2},h}^{k,r}+  10\cdot(1-  v^{k,r}_{n_2,f_2,h})  \hspace*{1cm}\forall k\in K, \forall r\in R_{k}, \forall h\in H_{k},\\
\hspace*{7.5cm} \forall n_{1}, n_2\in N, \forall f_{1},f_{2}\in F.
\end{array}\label{19}
\end{equation}
\begin{equation}
\begin{array}{l}
\displaystyle y_{n,f,h}^{k,r}\leq v_{n,f,h}^{k,r}  \hspace*{5.3cm} \forall k\in K, \forall r\in R_{k}, \forall h\in H_{k},\\ 
\hspace*{7.5cm} \forall n\in N, \forall f\in F.
\end{array}\label{20}
\end{equation}
\vspace*{0.05cm}

Constraints~\eqref{17}--\eqref{19} jointly ensure that at most one of the two terms to the left of Constraint~\eqref{17} is nonzero. This means that different feeds or feeds from different customers cannot be mixed in the same compartment of a truck. Constraint~\eqref{20} represents the logical relationship between $y$ and $v$.\\
}

\item[6)] {Constraints describing the nature of the decision variables.
\begin{equation}
     x_{n_{1},n_{2}}^{k,r} \in \{0,1\} \hspace*{6cm}\forall k\in K, \forall r\in R_{k}, \forall n_{1}, n_{2}\in \bar{N}.\label{21}
\end{equation}
\begin{equation}
     y_{n,f,h}^{k,r} \in [0,1] \hspace*{4.2cm}\forall k\in K, \forall r\in R_{k}, \forall n\in N, \forall f\in F, \forall h\in H_{k}.\label{22}
\end{equation} 
\begin{equation}
     z_{n,f}^{k,r} \in [0,1] \hspace*{5.7cm}\forall k\in K, \forall r\in R_{k}, \forall n\in \bar{N}, \forall f\in F.\label{23}
\end{equation}
\begin{equation}
    v_{n,f,h}^{k,r} \in \{0,1\} \hspace*{4cm}\forall k\in K, \forall r\in R_{k}, \forall n\in N, \forall f\in F, \forall h\in H_{k}.\label{24}
\end{equation}
\begin{equation}
    u_{n_{1},n_{2},f_{1},f_{2},h}^{k,r} \in [0,1] \hspace*{1.9cm}\forall k\in K, \forall r\in R_{k}, \forall n_{1}, n_{2}\in N, \forall f_{1},f_{2}\in F, \forall h\in H_{k}.\label{25}
\end{equation}
\vspace*{0.05cm}

Constraints~\eqref{21} and \eqref{24} specify that variables $x$ and $v$ are binary, and Constraints~\eqref{22}, \eqref{23}, and \eqref{25} indicate that variables $y$, $z$, and $u$ take values in $[0,1]$, respectively.

}
\end{enumerate}

\subsubsection{Objective function}\label{sec:obj}
We model the MC-VRPSD as a bi-objective optimization problem in the following way:

\begin{equation}
    \hspace{-4cm}
	\displaystyle \max\sum_{k\in K}\sum_{r\in R_{k}}\sum_{h\in H_{k}}\sum_{f\in F}\sum_{n\in N}c_{h}^{k}y_{n,f,h}^{k,r}\, ,\hspace*{0cm}
\label{eq:obj1}
\end{equation}

\begin{equation}
	\displaystyle \min
	\displaystyle \sum_{k\in K}\sum_{r\in R_{k}}
	\left(
	\sum_{n_{1}\in \bar{N}}
	\sum_{n_{2}\in \bar{N}} 
	d_{n_{1},n_{2}}x_{n_{1},n_{2}}^{k,r}  
	+ \sum_{n\in \hat{N}}
	\sum_{f\in \hat{F}_n}
	2 d_{0,n}  z_{n,f}^{k,r}
	\right).
\label{eq:obj2}
\end{equation}
\vspace*{0.05cm}

The objective function in~\eqref{eq:obj1} maximizes the total feed delivered during the current working day, while the objective function in~\eqref{eq:obj2} minimizes the transportation costs incurred by the cooperative. As we will see in Section~\ref{sec:results}, we will use the weighted sum method to address this multi-objective problem.

\section{Two-stage heuristic algorithm}\label{sec:heuristic}

Constructive algorithms and tabu search algorithms have proven to be effective tools for solving many NP-hard combinatorial optimization problems. Constructive algorithms create good initial solutions, while tabu search algorithms enhance these initial solutions 
by transiting from one solution to another through a number of iterations. This process involves making moves within a defined neighborhood structure and aims to avoid local optima and re-visiting previously explored solutions. We propose a two-stage heuristic algorithm for the MC-VRPSD presented in this paper, where the first phase consists of a constructive heuristic to generate a initial solution, and the second phase employs a tabu search algorithm. We integrate this two methods, resulting in an iterated tabu search (ITS) metaheuristic. The following subsections describe these procedures.

\subsection{Constructive heuristic}\label{sec:constructive}

In this section, we present the constructive heuristic algorithm designed to generate an initial solution. This algorithm is a generalization of the Clarke and Wright algorithm tailored to the specific characteristics of our stochastic problem, which includes compartmentalized vehicles and restrictions on accessibility, the urgency of certain orders, and the mixing of different products. The general idea is to select a customer and consider the route connecting the depot to that customer and back. In addition, a truck is selected to serve that customer, a compartment is assigned, and the compartment is loaded based on the customer's demand. Then, customers are inserted into the route using the criterion of maximum savings.
It is important to note that a customer may have a very large demand, which means that if only one compartment is assigned and the order is urgent, the truck will need to return to the cooperative headquarters after arriving at the customer's location to complete the order. For this reason, we will introduce the so-called customer replicas, which we explain below.\\

\noindent \textbf{Creation of customer replicas.} Given an urgent customer $n\in N$ and a feed $f\in F$, let $p_{n,f}=\mathbb{P}[\textbf{U}_{n,f}=0]$ represent the probability that the customer's order for that feed is urgent. Additionally, let $q_{n,f}$ denote the quantile of order $p_{n,f}$ of the random variable $\textbf{O}_{n,f}$, which corresponds to the stochastic demand of that customer for that feed. Thus, $q_{n,f}$ is a value for the demand that increases with the likelihood of the order being urgent, and represents a quantity of demand that is convenient to serve. 
Let $c=\max\{c^k_h\,:\,k\in K, h\in H_k\}$ be the maximum size among the compartments of all trucks, and define $r_{n,f}=\left\lceil \frac{q_{n,f}}{c} \right\rceil$.\footnote{$\left\lceil x \right\rceil$ denotes the smallest integer greater than or equal to $x$.} The value $r_{n,f}$ indicates the number of compartments needed to meet the demand $q_{n,f}$. 
For the convenience of algorithm design, customer $n$ will be replicated into $r_{n,f}$ customers (where $r_{n,f}\geq 1$) for each feed $f$. The distance between a replica of $n$ and a replica of another customer $n'$ will be equal to the distance between $n$ and $n'$. Furthermore, the distance between two replicas of the same customer will be considered as zero. On the other hand, the mean and standard deviation of the demand of each replica of customer $n$ for feed $f$ will be taken as the mean and the standard deviation of the demand of $n$, respectively, divided by $r_{n,f}$. Consequently, we consider a new problem with a total of $\sum_{(n,f)\in N\times F}r_{n,f}$ customers. Let $N^{*}$ denote the set comprising these customers. \\

\noindent \textbf{Selection of the first customer of a route.} 
After creating the customer replicas using the aforementioned procedure, the next step in the algorithm is to select the first customer to initiate the route. We prioritize urgent customers with large demands who are located farther from the depot. More specifically, when determining the first customer to be served, we take into account the distance to and from the depot (which is needed if we arrive at an urgent customer with insufficient load), the probability that the customer's order is urgent, and the probability that their demand exceeds the capacity of the assigned compartment. To assess the latter probability, we consider the probability that the demand exceeds the capacity of the compartment with the least capacity of all trucks, minus the probability that the demand exceeds the capacity of the compartment with the most capacity of all trucks. This quantity represents the probability that a customer's demand falls between the capacities of the smallest and the largest compartments. 
Formally, if $n\in N^*$ and $f\in F$, we define:

\[\displaystyle dif_{n,f}=2 d_{0,n}\left(\max_{k\in K, h\in H_{k}}\{\mathbb{P}[\textbf{O}_{n,f}>c^h_k]\}-\min_{k\in K,h\in H_{k}}\{\mathbb{P}[\textbf{O}_{n,f}>c^h_k]\}\right)\mathbb{P}[\textbf{U}_{n,f}=0].\]
We will select the customer $n$ and the feed $f\in F$, such that: $$\displaystyle{n\in \{n^*\in N^*\,:\,dif_{n^*,f}=\max_{(n',f')\in N^*\times F}dif_{n',f'}\}}.$$

Once $n\in N^*$ and $f\in F$ have been selected, the next step involves creating a route from the depot to the customer $n$ and back to the depot. For this, we select the vehicle $k_{n,f}$ and the compartment $h_{n,f}$ that are the most suitable because they provide more capacity. That is, these are chosen such that:
$$\mathbb{P}[\textbf{O}_{n,f}>c_{h_{n,f}}^{k_{n,f}}]=\min \{\mathbb{P}[\textbf{O}_{n,f}>c_{h}^{k}]\,:\, k\in K,\, h\in H_k\}.$$

Regarding the loading process, the compartment $h_{n,f}\in H_{k_{n,f}}$ will be filled with $q_{n,f}$, as long as $q_{n,f}\leq c_{h_{n,f}}^{k_{n,f}}$. If $q_{n,f}>c_{h_{n,f}}^{k_{n,f}}$, the load will be set to $c_{h_{n,f}}^{k_{n,f}}$, meaning the compartment will be filled to its maximum capacity. This approach aligns with the two objectives of our optimization problem: maximizing the amount of feed delivered while minimizing the total transportation costs.\\

\noindent \textbf{Incorporation of customers to a route.} In the next stage, another replica of customer $n$ is selected. If no such replica exits, we calculate, for each $n'\in N^* \backslash \{n\}$, the potential savings from inserting customer $n'$ (after $n$) to meet their demand for feed $f'$ on the current route:
\begin{equation}
\displaystyle S^{h_{n,f},k_{n,f}}_{n,n',f'}=d_{n,0}+d_{0,n'}-d_{n,n'}+\lambda I^{k_{n,f}, h_{n,f}}_{n',f'},
\label{saving}
\end{equation}
\noindent such that $I^{k_n, h_n}_{n',f'}$ is defined as:
\begin{equation}
\displaystyle 2 d_{0,n'}\left(\min_{k^{'}\in K\backslash\{k_{n,f}\}, h\in H_{k^{'}}}\{\mathbb{P}[\textbf{O}_{n',f'}>c_h^{k'}]\}-\min_{h\in H_{k_{n,f}}\backslash \{h_{n,f}\}}\{\mathbb{P}[\textbf{O}_{n',f'}>c_h^{k_{n,f}}]\}\right)\mathbb{P}[\textbf{U}_{n',f'}=0].
\label{saving2}
\end{equation}

In Equation~\eqref{saving2}, we evaluate the difference between the expected distance traveled (to and from the depot) when serving customer $n'$ with feed $f'$ using the largest available compartment of any truck, other than $k_{n,f}$ and the expected distance when using the largest free compartment of $k_{n,f}$, if the order is urgent. 
This difference represents the expected saving of serving $n'$ on the same route as $n$ instead of on a separate route. It accounts for the potential need to make an additional round trip to the depot if a smaller-than-needed load was delivered to $n'$ on the first visit. The value is positive when the compartment best suited to serve $n'$ is on truck $k_{n,f}$.
In Equation~\eqref{saving}, this amount is weighted by a factor $\lambda \in [0,+\infty)$, a parameter initialized at the start of the algorithm. 

\begin{figure}[ht!]
     \centering
     \begin{subfigure}[b]{0.49\textwidth}
         \centering
         \includegraphics[width=0.7\textwidth]{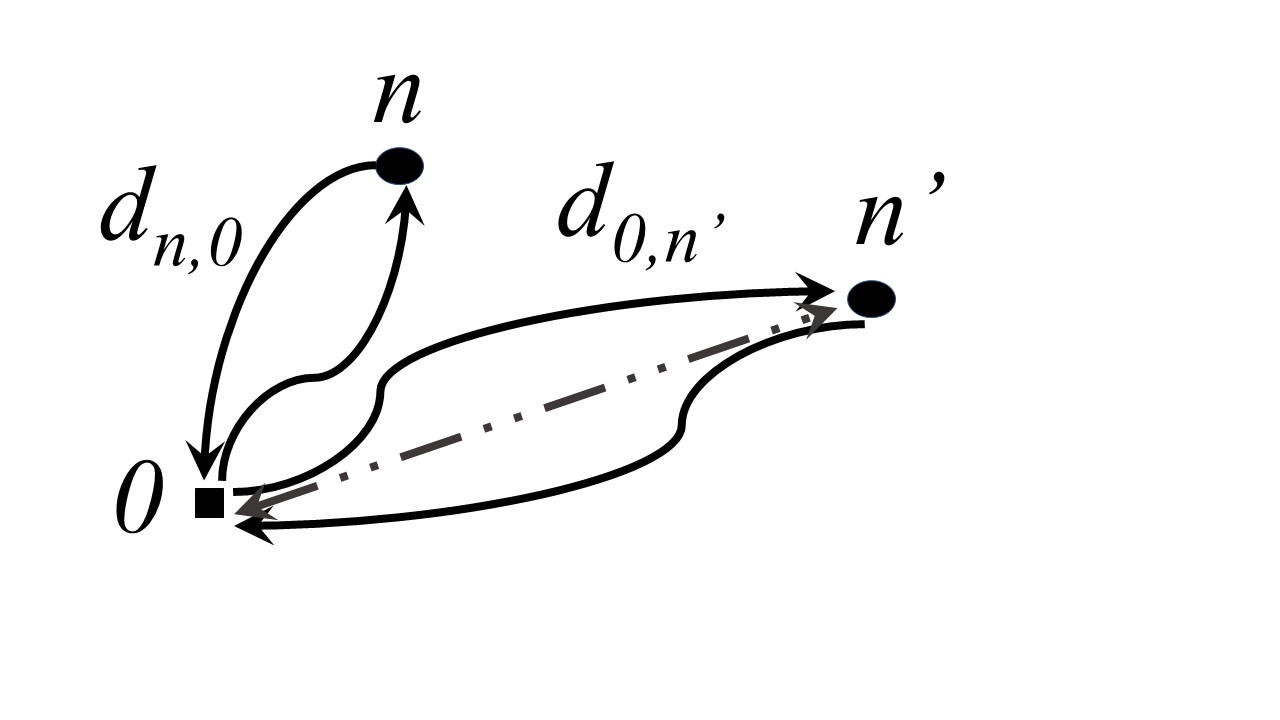}
         \caption{Before inserting $n'$ to the route of $n$.}
         \label{fig:cw1}
     \end{subfigure}
     \hfill
     \begin{subfigure}[b]{0.49\textwidth}
         \centering
         \includegraphics[width=0.7\textwidth]{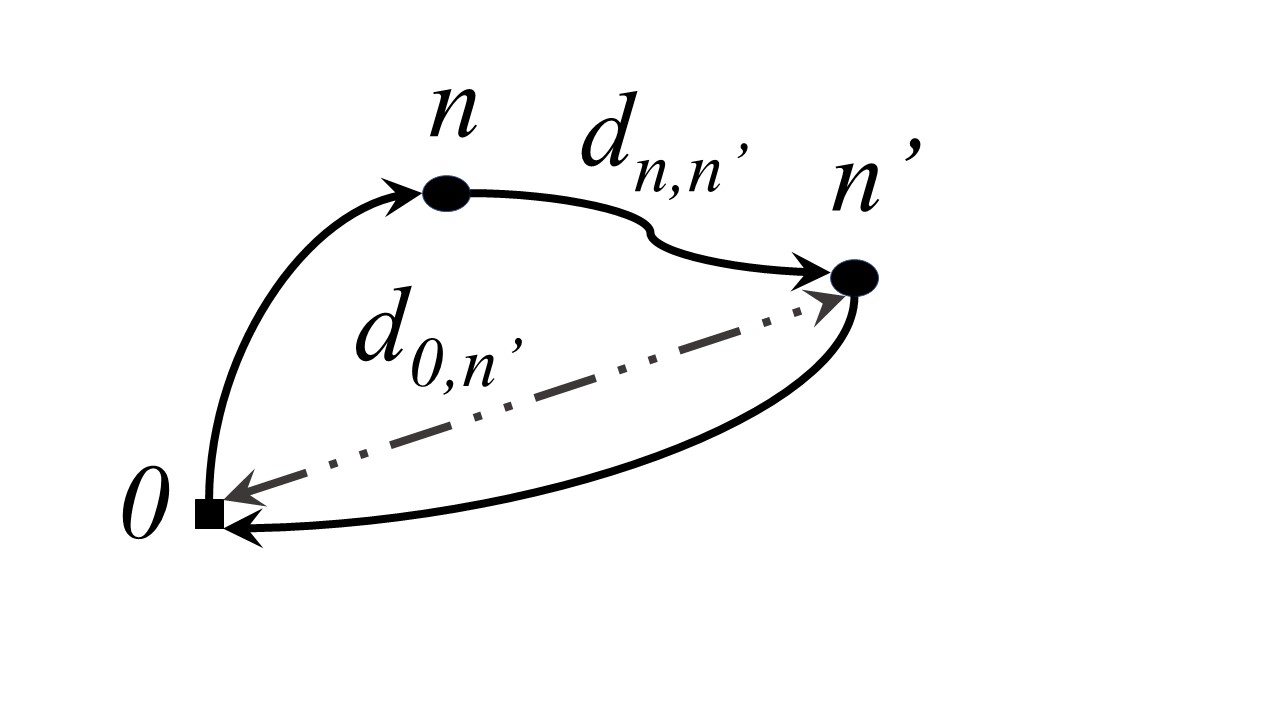}
         \caption{After inserting $n'$ to the route of $n$.}
         \label{fig:cw2}
     \end{subfigure}
             \caption{Calculation of $S^{h_{n,f},k_{n,f}}_{n,n',f'}$ for customers $n, n'\in N^{*}$.}
        \label{fig:cw}
\end{figure}

In summary, when evaluating the savings from inserting a new customer to a route, we consider both the fixed distances the truck travels and the expected distance required to return to the depot, while accounting for product separation constraints and the urgency of orders. Figure~\ref{fig:cw} illustrates the savings considered in this selection process.

Next, the customer that maximizes the savings defined by Equation~\eqref{saving} will be selected. Once a new customer $n'$ has been chosen using one of the methods described (either by selecting another replica of a customer or by identifying the customer that maximizes savings), it will be inserted to the route to be visited after $n$. This new customer is served by allocating the new load to one of the empty compartments, following the same procedure used for the first customer.

To finish this stage, we have to check the capacity restrictions and the maximum distance that each vehicle can travel. If both conditions are met, customer $n'$ is confirmed to be added to the route of truck $k_{n,f}$. If not, this stage is restarted by removing $n'$ from the set of candidates for insertion into the route.

The process of inserting customers into the route is repeated while the truck has free compartments. It is important to note that when inserting a new customer after having two or more customers per route (to broaden the search for the best insertion), we consider the savings $S_{n,n',f'}^{h_{n,f},k_{n,f}}$, where $n$ represents the last customer visited before returning to the depot on the current route, and the savings $S_{n',n,f'}^{h_{n,f},k_{n,f}}=d_{n',0}+d_{0,n}-d_{n',n}+\lambda I^{k_{n,f}, h_{n,f}}_{n',f'}$, where $n$ represents the first customer visited after leaving the depot on the current route. In both cases, $n'$ refers to any customer in $N^*$ who has not yet been assigned to a route and is a candidate for insertion either after $n$ or before $n$, respectively. Once the route of a truck has been completely designed, the other trucks will be assigned routes in a similar way.

Algorithm~\ref{ps:psoalg} schematically shows the procedure for planning a route for each truck.

\begin{algorithm}[ht!]
\small
	\caption{Savings-based construction heuristic for the MC-VRPSD.}\label{ps:psoalg}
\begin{algorithmic}[1]
\State{make replicas of customers taking into account demands; set $\lambda$, $K_0=K$}
\Comment{\textcolor{gray}{Initialization.}} \label{lin:ini}
\While{there are unattended customers}
\State{select customer $n\in N^{*}$, $f\in F$, who maximize $dif_{n,f}$}
\Comment{\textcolor{gray}{First customer of a route.}} \label{lin:first}
\State{create the route $r:=(0,n,0)$ for a vehicle $k_{n,f}\in K$ and allocate the load}
\If{there are unattended customers}
\Comment{\textcolor{gray}{Route completion.}} \label{lin:ini3}
\State{select a new replica of $n$ or $n'\in N^{*}$, $f'\in F$, that maximize $S^{h_n,k_n}_{n,n',f'}$ or $S^{h_n,k_n}_{n',n,f'}$}
\If{capacity and distance constraints are satisfied}
\State{insert $n'$ in the route and return to~\ref{lin:ini3}}
\State{$n \gets n'$}
\Else
\State{discard $n'$ and go back to~\ref{lin:ini3}}
\EndIf
\State{end~\ref{lin:ini3} when the vehicle is full or there are no feasible insertions}
\EndIf
\State{$K_0 \gets K_0\backslash \{k_{n,f}\}$}
\EndWhile
\end{algorithmic}
\end{algorithm}

\normalsize

\noindent \textbf{Example of how the constructive heuristic works.}  
In this paragraph, we illustrate the performance of our algorithm using a small fictional example. We consider a problem with three customers (1, 2, and 3) and a depot (0). The distances, in minutes, between the depot and each customer and between the customers themselves, are presented in Table~\ref{sample-table}.

\begin{table}[ht!]
\centering
\caption{Travel distances between pairs of nodes in our fictitious example.}
\label{sample-table}       
\begin{tabular}{lcccc}
\toprule
Node & 0 & 1 & 2 & 3   \\
\midrule
$0$  & -- & $28$ & $69$ & $64$  \\  
		$1$  & $28$ & -- & $67$ & 62  \\ 
		$2$  & $69$ & $67$ & -- & $7$  \\  
		$3$  & $64$ & $62$ & $7$ & --  \\   
			\botrule
\end{tabular}
\end{table}

The average order sizes of the different customers are $3.30$, $2.95$, and $3$ tonnes, respectively.\footnote{For simplicity, since each customer only orders one type of feed, we will omit the subscript $f$ indicating the feed type from the parameters $p_{n,f}$, $q_{n,f}$, $r_{n,f}$, and $dif_{n,f}$ throughout the example, using $p_{n}$, $q_{n}$, $r_{n}$, and $dif_{n}$ instead, respectively.} We assume that the orders are random variables following a normal distribution with a standard deviation of $0.5$. Additionally, we consider that all customers have urgent orders for only one type of feed, with $\mathbb{P}(\textbf{U}_{n,f}=0)=0.95$ for each $n\in \{1,2,3\}$ and $f\in F$.
Furthermore, we take one truck with five compartments, which can hold up to $3$, $3.7$, $3.8$, $3.7$, and $3$ tonnes of feed, respectively. The total carrying capacity of the truck is limited to $11.8$ tonnes, and the truck can deliver food to all customers.

\begin{figure}[ht!]
     \centering
              \includegraphics[width=0.85\textwidth]{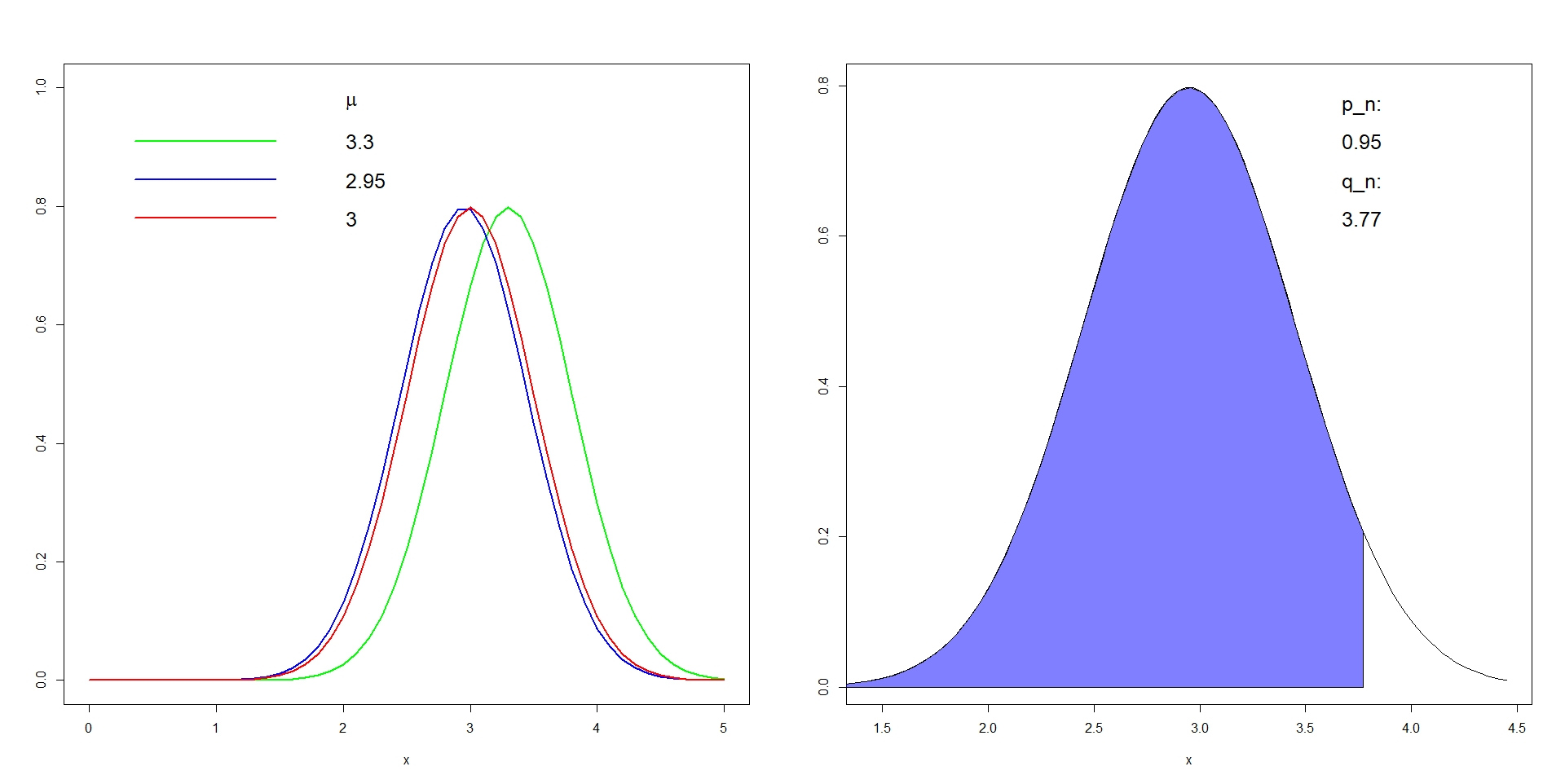}
         \caption{Stochastic demands of three customers (left) and $0.95$-th quantile of customer~2's demand (right): $q_n=3.77$ will be the load allocated to customer~2 in the largest compartment of the truck.}
         \label{fig:cwbis}
\end{figure}

In this case, $p_n=0.95$ for each $n\in \{1,2,3\}$. Considering the statistical distribution of the demands (see Figure~\ref{fig:cwbis}, left-hand side), we obtain $q_1=4.12$, $q_2=3.77$ (see Figure~\ref{fig:cwbis}, right-hand side) and $q_3=3.82$. Given that $c=3.8$, we have $r_1=r_3=2$ and $r_2=1$. Consequently, customer~1 is replicated in two customers with a demand following a normal distribution with an average of $3.30/2=1.65$ and a standard deviation of $0.5/2=0.25$, and customer~3 is replicated in two customers with a demand following a normal distribution with an average of $3/2=1.5$ and a standard deviation of $0.5/2=0.25$. Note that customer~2 has the largest average demand in the new set $N^*$ of five customers and is also the furthest from the depot. Consistently, we obtain
$dif_2=\max_{n'\in N^*}dif_{n'}$. Therefore, a route is created that starts from the depot, goes to customer~2, and returns to the depot. In addition, the compartment with the largest capacity of the vehicle (compartment number~3) is selected and filled with a load of $q_2=3.77$.

\begin{figure}[ht!]
     \centering
     \begin{subfigure}[b]{0.49\textwidth}
         \centering
         \includegraphics[width=0.7\textwidth]{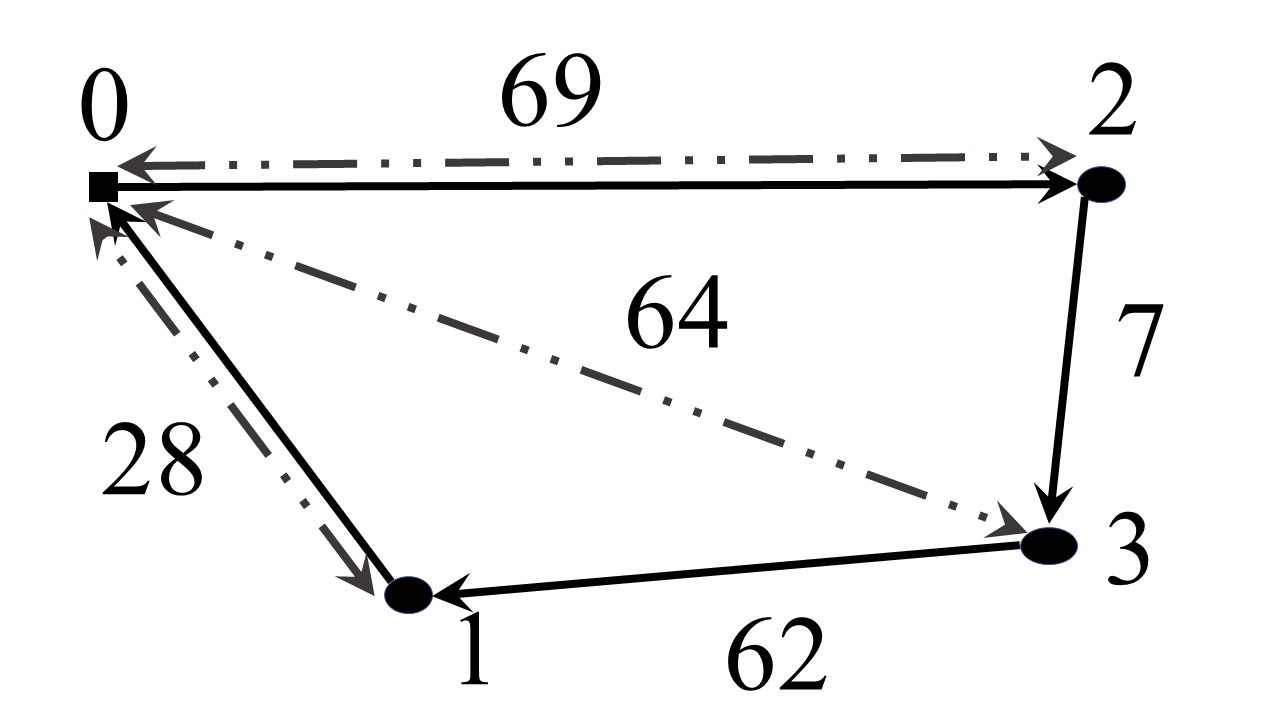}
         \caption{Route (order of customer visits and distances).}
         \label{fig:cw3}
     \end{subfigure}
     \hfill
     \begin{subfigure}[b]{0.49\textwidth}
         \centering
         \includegraphics[width=0.7\textwidth]{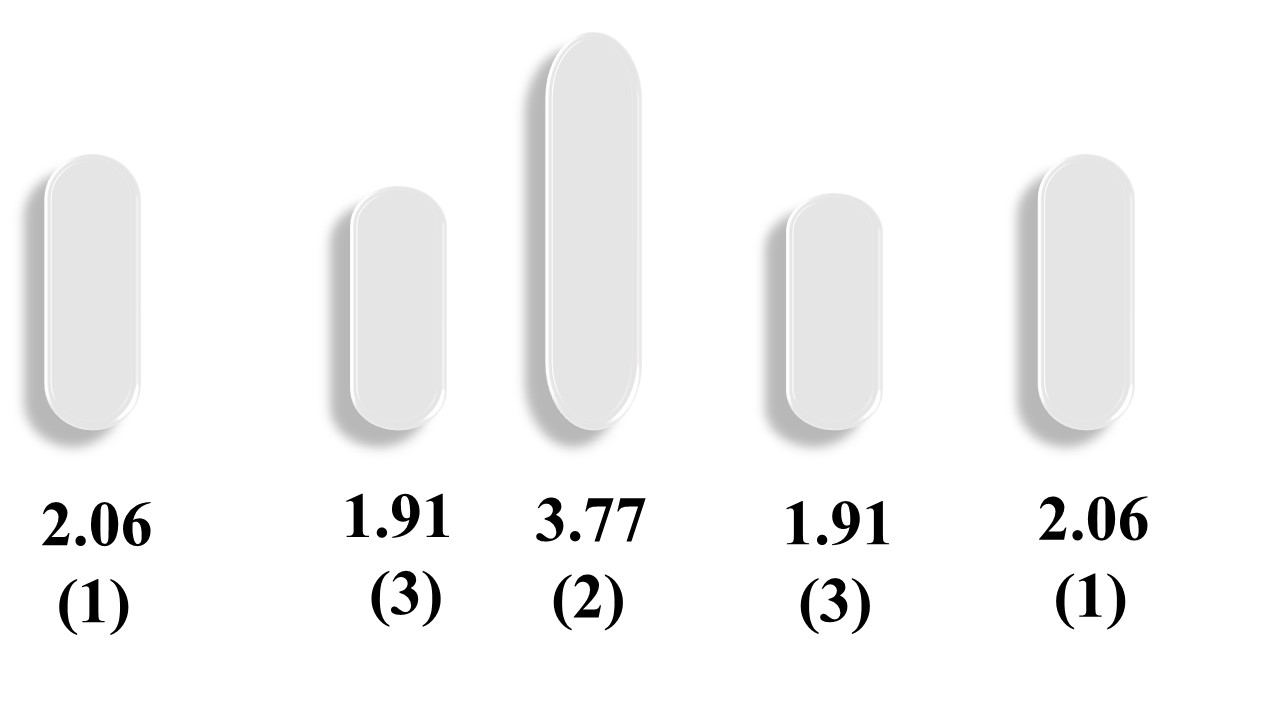}
         \caption{Load, in tonnes, and customer allocation (in brackets) of the compartments.}
         \label{fig:cw4}
     \end{subfigure}
             \caption{Solution obtained with the constructive heuristic.}
        \label{fig:cw34}
\end{figure}

Now, if $n^*$ is a replica of customer~1, considering they have a normal demand with an average of $1.65$ and a standard deviation of $0.25$, we have that $q_{n^*}=2.06$. Similarly, if $n^*$ is a replica of customer~3, considering they have a normal demand with an average of $1.5$ and a standard deviation of $0.25$, we have that $q_{n^*}=1.91$. Using Equation~\eqref{saving} with $\lambda=1$, we first add customer~3 to the route, assigning them the two free compartments with the largest capacity. Subsequently, customer~1 is incorporated and allocated the two remaining compartments. Figure~\ref{fig:cw3} illustrates the resulting route, while Figure~\ref{fig:cw4} displays the allocation of compartments to customers and the quantity loaded into each compartment. As a result, the truck departs from the depot with a total load of $11.71$ tonnes. The total distance traveled will be $182.3=166+16.3$, where $166$ represents the fixed distance to be covered, and $16.3$ is the expected distance due to the probability of reaching a customer with a quantity for them which is less than their urgent needs.

\subsection{Iterated tabu search metaheuristic}\label{sec:its}

In this section, we introduce an iterated tabu search (ITS) designed to improve the solution provided by the constructive algorithm presented in the previous section. Our metaheuristic is an adaptation of the approach developed by \cite{Sil17}, one of the best recent metaheuristics for solving MC-VRPs. Specifically, we implement a tabu search embedded in an iterated local search framework. In essence, our ITS consists of an initialization step followed by a tabu search and a perturbation phase that modifies the solution generated by the tabu search, and which are executed iteratively.
The tabu search (\citealp{Glo89, Glo90}) is a local search technique that manages to improve its performance through a memory structure. This structure ensures that once a potential solution is found, it is
not revisited in the immediate iterations of the algorithm. As a local search method, it makes use of a neighborhood structure that allows moving from one solution to another, until a specified stopping criterion is met. Different tabu search variants exist in the literature, and we have adapted the version proposed by \cite{Os93} for the VRP, as their approach has been shown to reduce both computational time and travel distances. Regarding the perturbation phase, we design and implement a variety of destroy and repair operators based on the proposals of \cite{Ali18} for the multi-depot MC-VRP. Destroy operators remove part of the solution, while the repair operators reconstruct the solution in a different way.\\

\noindent \textbf{ITS initialization.} Our iterated tabu search algorithm begins by selecting five parameters: the number of routes to be created, $\sigma$; the number of customers to be destroyed, $s$; the maximum number of exchanges allowed between two routes, $\kappa$; the maximum number of iterations for the tabu search, $m$; and the maximum number of iterations for the perturbation phase, $\iota$. Additionally, an initial solution is required, which is provided by the constructive algorithm presented in the previous section.\\

\begin{figure}[ht!]
     \centering
     \begin{subfigure}[b]{0.49\textwidth}
         \centering
         \includegraphics[width=0.7\textwidth]{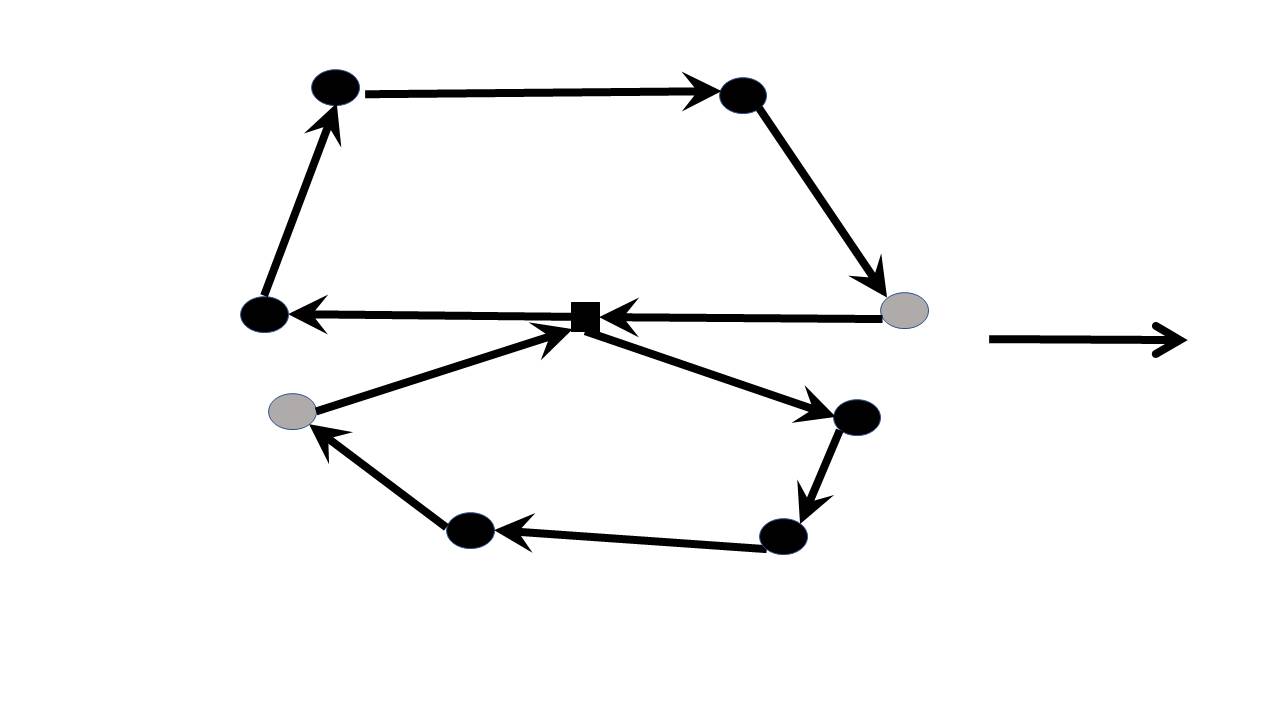}
         \caption{}
         \label{fig:k1}
     \end{subfigure}
     \hfill
     \begin{subfigure}[b]{0.49\textwidth}
         \centering
         \includegraphics[width=0.7\textwidth]{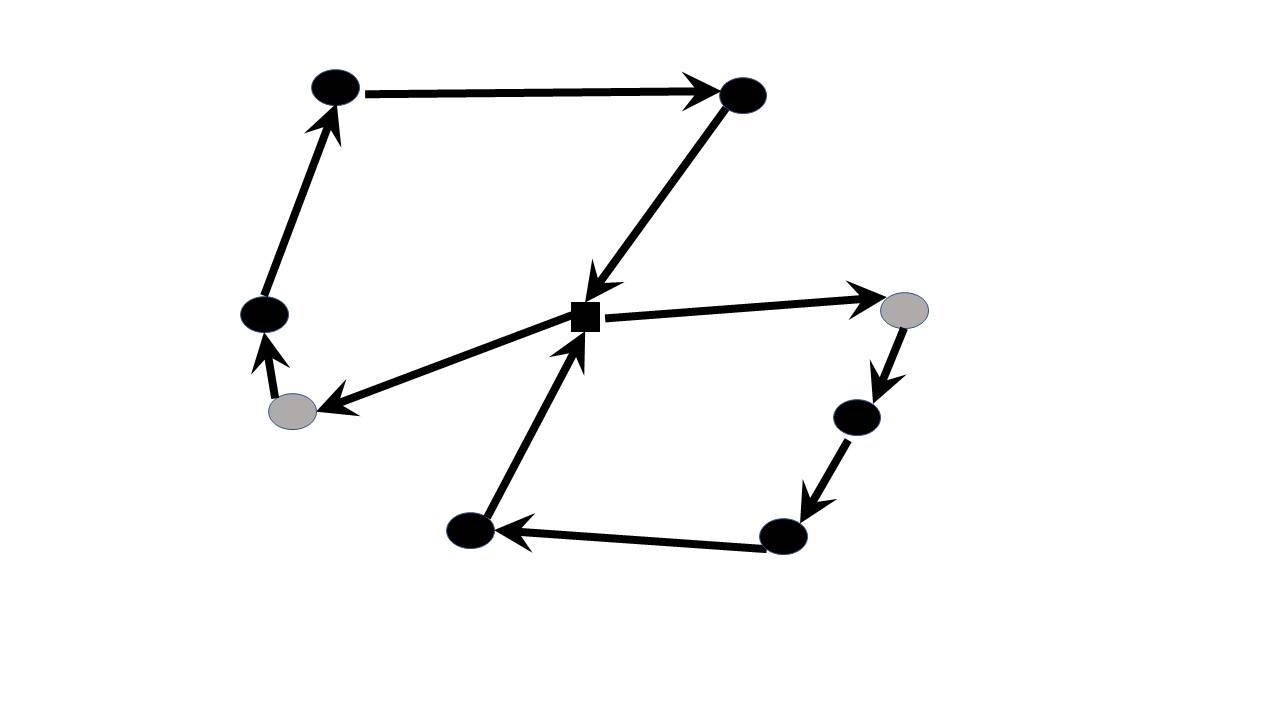}
         \caption{}
         \label{fig:k1bis}
     \end{subfigure}
             \caption{Example of a 1-exchange between two routes, generating two neighboring solutions.}
        \label{fig:k11}
\end{figure}

\noindent \textbf{Neighborhood structure.} The neighborhood structure is based on the so-called $\kappa$-exchanges. More specifically, consider two different routes, $r_1$ and $r_2$, which might be assigned to either two different trucks or to separate routes of the same truck in a given solution. A $\kappa$-exchange consists in transferring up to $\kappa$ customers from $r_1$ to $r_2$, and up to $\kappa$ customers from $r_2$ to $r_1$. The number of customers transferred between the routes does not need to be equal, but both transfers must involve no more than $\kappa$ customers. Additionally, one route may transfer customers without receiving any in return. Figure~\ref{fig:k11} exemplifies this neighborhood structure. After performing a $\kappa$-exchange, a random number of customer rearrangements within each route is evaluated  in order to reduce the distance traveled. Furthermore, the use of the vehicle's compartments is re-optimized by prioritizing urgent customers first, aiming to minimize the expected distance traveled. Afterward, the remaining customers are  assigned to compartments based on their demand, with larger-capacity compartments allocated to those with higher demand. Also, note that, if two replicas of the same initial customer are present on the same route, it is possible to move only one to another route. For practical purposes, replicas corresponding to the same customer are treated as different customers located at the same site.

\noindent \textbf{Tabu search.} In each iteration, we have two options for selecting the different routes to perform exchanges. The first option involves arbitrarily choosing $\sigma$ routes from the current solution. The second option consists of randomly choosing one route and identifying the nearest customer not included in that route for each customer on it. This process allows us to select $\sigma-1$ additional routes that contain these nearest customers. A neighboring solution is then obtained by performing a $\kappa$-exchange between two of the selected $\sigma$ routes. Next, all feasible neighboring solutions that satisfy the problem constraints are considered. Finally, the solution with the smallest total distance is selected, updating the current solution, and the search progresses to the next iteration.

At the end of each iteration, information about the exchange performed is stored in the tabu list. Specifically, each customer that has been transferred from one route to another is recorded. These customers are restricted from moving to any route, whether their original route or another, for three iterations.\footnote{This value is another parameter of the algorithm. Different values were tested and, finally, the value 3 was selected.} 
It should also be noted that an iteration may conclude with a current solution that does not yield a total distance lower than those achieved previously. Allowing this type of iteration helps prevent the algorithm from getting stuck in local optima; however, a maximum of three such iterations can be performed. Additionally, there might be cases where no feasible exchanges can be made. In this situation, the best solution found up to that point is taken, and the algorithm continues until the specified number of iterations is reached. 

Finally, a $2$-opt local search algorithm, based on the GENI method proposed by \cite{Gen92}, is applied to improve the routes following the tabu search.\\

\noindent \textbf{Perturbation phase.} In the second step of the algorithm, a destruction and repair method is applied. For the destruction method, three different ways of selecting customers are proposed:

\begin{enumerate}

    \item $s$ customers are arbitrarily chosen from the current solution (random removal). 

    \item One customer is arbitrarily chosen, and the $s-1$ nearest customers to that customer are also selected (Shaw removal). 

    \item The $s$ customers with the highest individual on-route cost $C^{k,r}_{n,f,h}$ are chosen, where:
    \[C^{k,r}_{n,f,h}=d_{n_{1},n}+d_{n,n_{2}}+2d_{0,n}P^{k,r}_{n,f,h}\]
    with
    \begin{equation*}
    P_{n,f,h}^{k,r}:=  \left\lbrace
    \begin{array}{lcl}
    0 & \; & \text{if } \mathbb{P}\left[ \textbf{U}_{n,f}=0\right]\leq \beta,\\
    \mathbb{P}\left[c_{h}^{k}y_{n,f,h}^{k,r}\leq \textbf{O}_{n,f}\right] & \; & \text{otherwise},
    \end{array}
    \right.\hspace*{4.5cm}
    \end{equation*}
    for all $n\in N,\ k\in K$, $r\in R_k,\   h\in H_{k},$ and $ f\in F$; where $n_{1}$ and $n_{2}$ are the predecessor and successor of customer $n$ in the route $r$, respectively, and $\beta=0.90$ (worst distance removal).

\end{enumerate}

Note that the second and third methods can be combined, such as by choosing the customer with the highest individual on-route cost and the $s-1$ nearest customers to that customer.
 
In the repair stage, for each customer selected in the destruction operation, we calculate the difference between the best and the second best objective values obtained by inserting that customer into the routes of the destroyed solution, exploring all possible insertion points. We then insert the customer with the highest difference into the position that results in the lowest cost, followed by a re-optimization of the compartment utilization. This process is repeatead with the updated solution until all customers have been inserted. 

Finally, after the perturbation phase, a $2$-opt local search algorithm is applied to improve the routes, using the GENI method as described by \cite{Gen92}.

Algorithm~\ref{ps:psoalg2} presents a schematic overview of the ITS described. 

\begin{algorithm}[ht!]
\small
	\caption{Iterated tabu search for the MC-VRPSD.}\label{ps:psoalg2}
\begin{algorithmic}[1]
\State{start with an initial solution provided by the constructive algorithm as the current solution and fix $\sigma$, $s$, $\kappa$, $m$, and $\iota$} 
\Comment{\textcolor{gray}{Initialization.}} \label{lin:ini2}
\State{set $iter=0$, $s.iter=0$, $p.iter=0$}
\While{$iter<m$}
\Comment{\textcolor{gray}{Tabu search.}} \label{lin:tabu}
\State{$iter=iter+1$}
\State{select $\sigma$ routes from the current solution}
\State{make all possible $\kappa$-exchanges}
\State{choose the feasible $\kappa$-exchange that minimizes the total distance and update the current solution}
\State{add exchanged customers to the tabu list and, if $iter\geq4$, remove the first entry from the tabu list}
\If{the new solution does not improve the previous solutions}
\State{$s.iter=s.iter+1$}
\If{$s.iter=3$}
\State{stop}
\EndIf
\EndIf
\EndWhile
\State{select the best solution obtained in~\ref{lin:tabu}}
\State{apply 2-opt local search algorithm}
\If{$p.iter<\iota$}
\Comment{\textcolor{gray}{Perturbation phase.}} \label{lin:desandrep}
\State{destroy current solution}
\State{repair destroyed solution}
\State{$p.iter=p.iter+1$}
\State{set $iter=0$, $s.iter=0$}
\State{apply 2-opt local search algorithm}
\State{return to~\ref{lin:tabu}}
\Else
\State{stop}
\EndIf
\end{algorithmic}
\end{algorithm}
\normalsize

\section{Results and discussion}\label{sec:results}

To evaluate the efficiency of our proposals, we analyzed the impact of the MC-VRPSD model and the two-phase heuristic using a modified version of a set of well-known benchmark problems from the literature, alongside a real-world problem. The proposed heuristic described in Section~\ref{sec:heuristic} was implemented in \texttt{R 4.0.2}. To validate its performance, we conducted a set of experiments. The mathematical model presented in Subsection~\ref{sec:formulation} was solved using the \texttt{Gurobi} 8.1.0 solver. The code was run on an Intel Core i5-10600K CPU 4.1GHz with 32 GB RAM.

The following subsections report the computational study. Subsection~\ref{sec:exact_model} provides a detailed examination of the exact solution of the proposed MC-VRPSD model, using first a numerical example and then a small example corresponding to the case study. This latter real-world application (described in Subsection~\ref{sec:exact_real}) was used to validate our heuristic, as reported in Subsection~\ref{sec:heuristic_small}, where we compare the quality of the results achieved with those from the exact method. 
To the best of our knowledge, there are currently no benchmark problems corresponding to our MC-VRPSD. Therefore, before presenting the computational results of our algorithm, Subsection~\ref{sec:new_instances} describes the generation of a set of new MC-VRPSD test problems, based on the 14 well-known MC-VRP cases introduced by \cite{elf08}. Finally, Subsection~\ref{sec:heuristic_results_instances} presents the solutions obtained using our proposed heuristic for both the instances created by \cite{elf08} for the MC-VRP and the datasets we generated for the MC-VRPSD. We compare the results for the MC-VRP with those available in the literature.

\subsection{Exact solving of the model}\label{sec:exact_model} 

The MC-VRPSD model is illustrated below through a numerical example and a real-world case. In the latter, we compare the obtained results with those that would arise in the absence of uncertainty. In addition, we show how the computation time significantly increases with a larger number of customers whose demands the cooperative aims to meet, as well as the trade-off between the gap and  computation time. For this analysis, we utilize the \texttt{AMPL} modeler and the \texttt{Gurobi} solver. For the sake of simplicity, this case does not account for any accessibility restrictions, and only one trip per vehicle is permitted.

\subsubsection{Numerical example}\label{sec:numerical_example}

We consider a problem with four customers, $N=\{1, 2, 3, 4\}$, each requiring a single type of feed. Customers~1 and 2 have urgent demands. There is a vehicle equipped with four compartments with capacities of $7$, $6$, $6$, and $6$ units, respectively. The demands for customers~1 and 2 are modeled as discrete random variables, identically distributed, taking values $5$, $6$, and $7$ units with probabilities $5/10$, $4/10$, and $1/10$, respectively. Customers~3 and 4 have fixed demands of $7$ units each. The distance between any pair of customers, as well as between the cooperative headquarters and each customer, is consistently set at $10$ units. Our objective function is formulated as a weighted average of the two model objectives, using a parameter $\omega \in [0,1]$ as follows:

$$
\begin{array}{l}
	\displaystyle \min \ \omega
	\displaystyle \sum_{k\in K}\sum_{r\in R_{k}}
	\left(
	\sum_{n_{1}\in \bar{N}}
	\sum_{n_{2}\in \bar{N}} 
	d_{n_{1},n_{2}}x_{n_{1},n_{2}}^{k,r}  
	+ \sum_{n\in \hat{N}}
	\sum_{f\in \hat{F}_n}
	2 d_{0,n}  z_{n,f}^{k,r}
	\right)\\
	\\
	\displaystyle \hspace*{0.75cm} -(1-\omega)\sum_{k\in K}\sum_{r\in R_{k}}\sum_{h\in H_{k}}\sum_{f\in F}\sum_{n\in N}c_{h}^{k}y_{n,f,h}^{k,r}. \hspace*{6cm}\\
\end{array}
$$

\begin{figure}[ht!]
     \centering
     \begin{subfigure}[b]{0.3\textwidth}
         \centering
         \includegraphics[width=\textwidth]{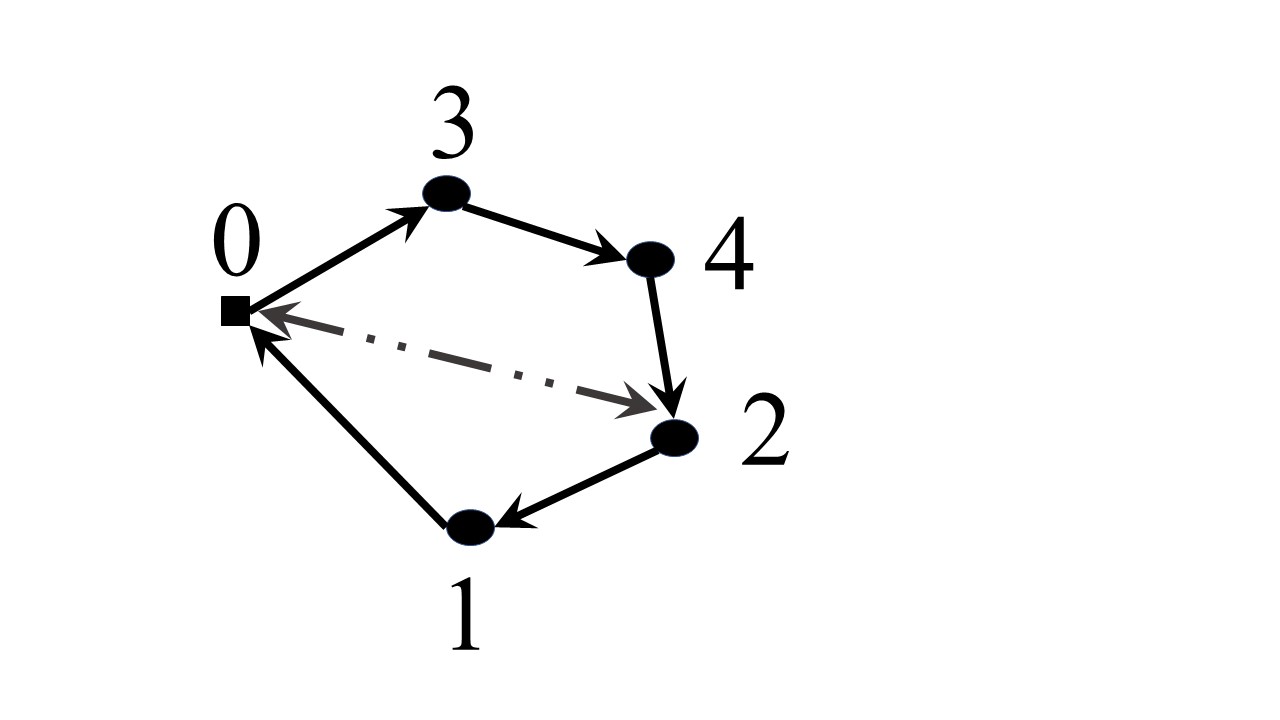}
         \caption{$\omega=0.2$.}
         \label{fig:w1}
     \end{subfigure}
     \hfill
     \begin{subfigure}[b]{0.3\textwidth}
         \centering
         \includegraphics[width=\textwidth]{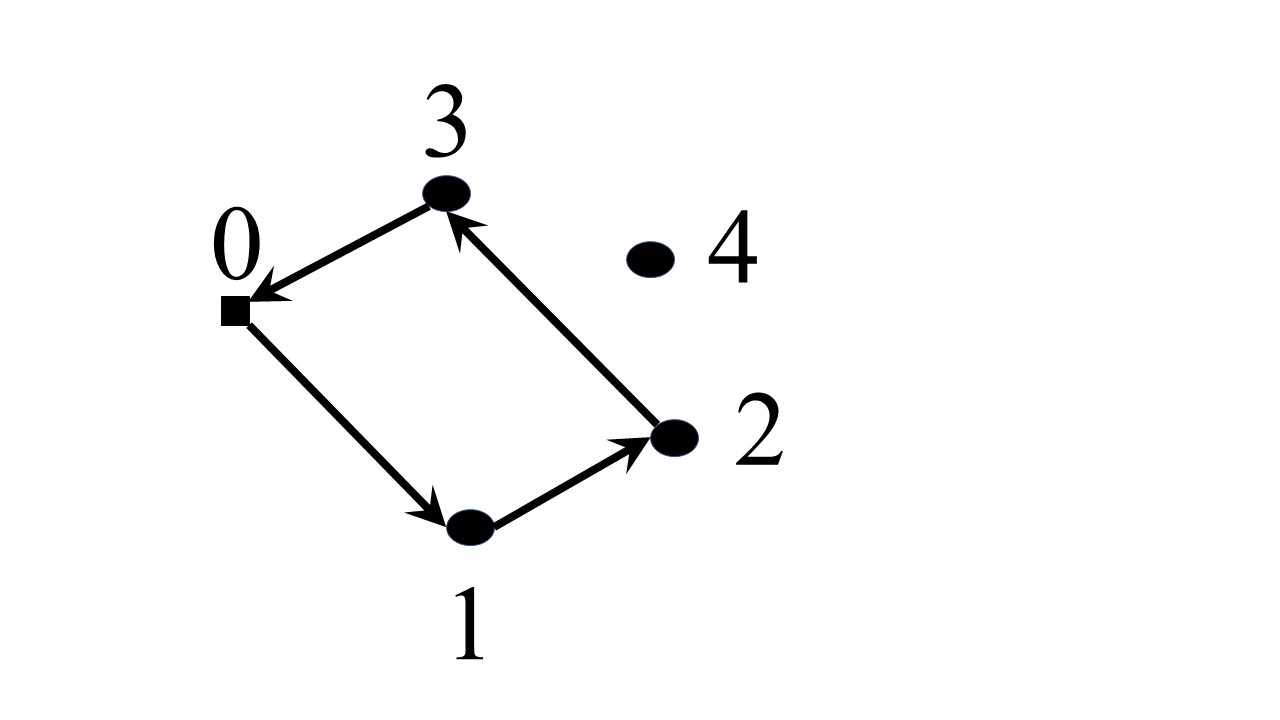}
         \caption{$\omega=0.3$.}
         \label{fig:w2}
     \end{subfigure}
     \hfill
     \begin{subfigure}[b]{0.3\textwidth}
         \centering
         \includegraphics[width=\textwidth]{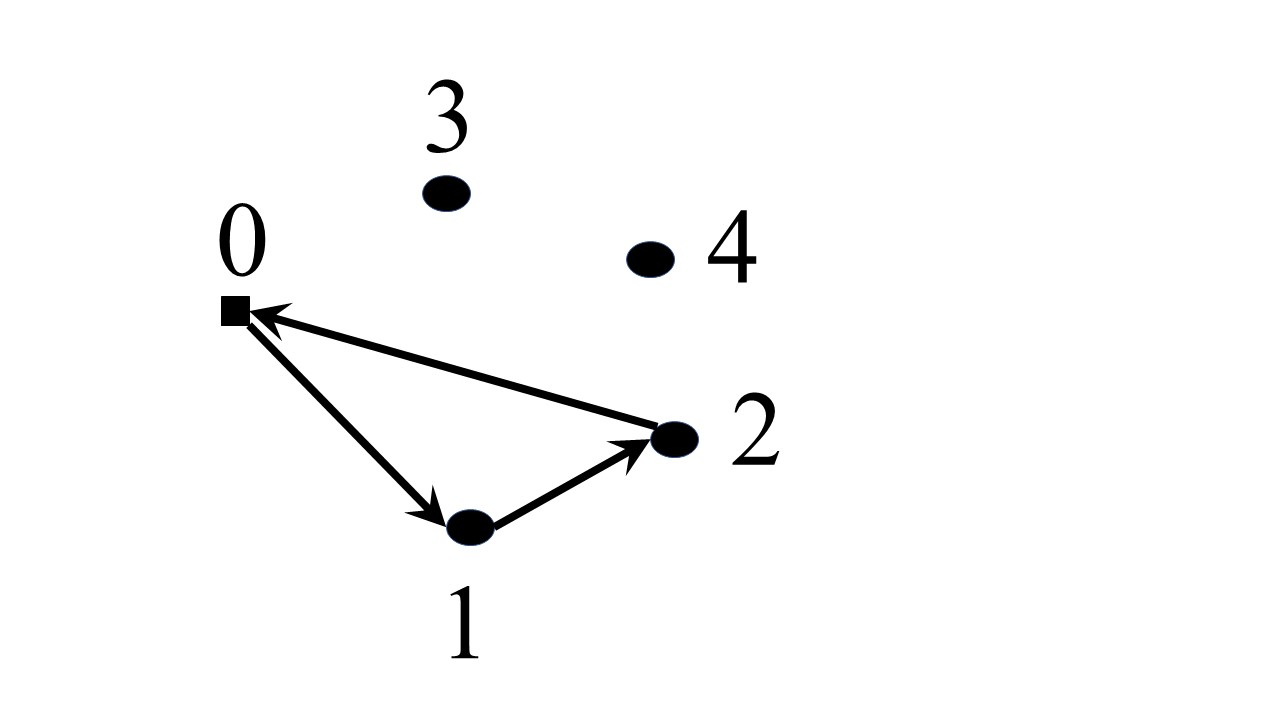}
         \caption{$\omega=0.8$.}
         \label{fig:w3}
     \end{subfigure}
        \caption{Optimal routes of the numerical example for varying parameters of $\omega$.}
        \label{fig:w}
\end{figure}

We present the solutions obtained for different values of the parameter $\omega$, which correspond to varying weights assigned to the two objectives. Initially, we set $\omega = 0.2$, prioritizing the maximization of the transported load. Next, we take $\omega = 0.3$, thereby increasing the importance given to the distance traveled. Finally, we take $\omega = 0.8$, significantly prioritizing the minimization of the distance traveled. Figure~\ref{fig:w} illustrates the optimal routes obtained for each scenario. Note that 0 represents the depot, and the dotted line indicates a potential round trip between a customer and the depot, facilitating the precise fulfillment of that customer's demand. 

Table~\ref{tab:example1} shows, for each scenario, the initial load allocated to each compartment of the truck, along with the corresponding customer for that load. Finally, Table~\ref{tab:example2} provides the value of the objective function, the expected total distance to be traveled, the transported load (excluding round trips between a customer and the depot), and the execution time of the solver. Together,  Figure~\ref{fig:w} and Tables~\ref{tab:example1} and \ref{tab:example2} illustrate how the different weights assigned in the objective function affect the results.

\begin{table}[!ht]
\centering
\caption{Compartment allocation in the numerical example.}\label{tab:example1}
\begin{tabular}{c|cccccc}
\toprule
              & \multicolumn{2}{c}{$\omega=0.2$} & \multicolumn{2}{c}{$\omega=0.3$} & \multicolumn{2}{c}{$\omega=0.8$}\\
Compartment  &   \multicolumn{2}{c}{\hrulefill} & \multicolumn{2}{c}{\hrulefill}   & \multicolumn{2}{c}{\hrulefill} \\
              &  Customer    & Load               & Customer     & Load               & Customer       &  Load   \\
\midrule 
   1          &     1       &       7            &      1      &         7          &    1          &   7   \\
      2          &     2       &       6            &      2      &         6          &    2          &   6   \\
         3          &     3       &       6            &      2      &         1          &    2          &   1   \\
            4          &     4       &       6            &      3      &         6          &    --          &   --   \\
            \botrule
            
\end{tabular}
\end{table}

\begin{table}[!ht]
\centering
\caption{Additional results in the numerical example.}\label{tab:example2}
\begin{tabular}{c|cccc}
\toprule
Weight ($\omega$) &  Objective    & Expected           & Total load  & Time   \\
                   &               &  overall distance  &             &  (seconds) \\
\midrule 
  0.2   &  -9.6    &   52  & 25 & 0.18 \\
  0.3   &  -2.0      &   40  & 20 & 0.05 \\
  0.8   &   21.2   &   30  & 14 & 0.06 \\ 
  \botrule
\end{tabular}
\end{table}

\subsubsection{Exact solving of a small real-world example}\label{sec:exact_real}

This section presents a study of the exact solution to the previously proposed MC-VRPSD model using the case of the agricultural cooperative that distributes feed to various farms (see Subsection~\ref{sec:description}). With real data available for this case study, we can construct an optimization scenario to asses our MC-VRPSD model. In particular, we have information on the distances between  different nodes (customers and the cooperative headquarters), as well as the demands of the customers and the capacities of the trucks. The drivers' working hours are 9 hours, and we estimate the average vehicle speed to be 60 km/h. Due to a lack of data on the service time for each customer or the time taken to load the trucks, we will assume these times are negligible. In this real-world instance of the model, we have selected 10 customers, served by one or two company trucks.
The trucks have five compartments that can contain up to $4000$, $3000$, $1700$, $4500$, and $3000$ kilograms (kgs), respectively. The maximum permissible load is $15,300$ kgs. Table~\ref{tab1} presents the matrix of distances (in minutes) between the nodes. In a stochastic scenario, the first five customers have equiprobable demands of $2990$, $3300$, $3500$; $5250$, $5500$, $6041$; $5560$, $5730$, $5959$; $2680$, $2951$, $3100$; and $4320$, $4490$, $4885$, respectively. In a deterministic case, their demands are $3300$, $6041$, $5959$, $2951$, and $4885$, all of which are urgent. The demands of the remaining five customers are $3003$, $3016$, $4478$, $5413$, and $3490$, respectively.

\begin{table}[ht!]
\centering
\caption{Matrix of distances (in minutes) between the nodes.\label{tab1}}
\begin{tabular}{c|ccccccccccc}
&\textbf{0}	& \textbf{1}	& \textbf{2}& \textbf{3}& \textbf{4}& \textbf{5}& \textbf{6}& \textbf{7}& \textbf{8}& \textbf{9}& \textbf{10}\\
\hline 
\textbf{0}	& 0     &     21   &       20   &       17      &    65      &    63&          60       &   19     &     22  &        24  &        60   \\     
\textbf{1}	&21      &     0    &       4    &       6  &        60  &        58&          55   &       15     &     18    &      20  &        55  \\
\textbf{2}	&20 &          4    &       0    &       4     &     59    &      56        &  53    &      13        &   8    &      12     &     53\\ 
\textbf{3}	&17   &         6  &          4  &          0   &        57 &          54        &   52  &         11   &        13   &        16 &          52\\
\textbf{4}	&65      &    60   &       59   &       57   &        0    &       3& 
7    &       66   &        69    &       71   &         6   \\     
\textbf{5}	&63   &       58  &        56  &        54  &         3     &      0&           4      &    64   &       66  &        69   &        3\\
\textbf{6}	&60 &         55   &       53   &       52  &         7    &       4&           0     &     61    &      64     &     66    &       2\\ 
\textbf{7}	&19   &       15    &      13    &      11   &       66   &       64&          61       &    0        &   3    &       5     &     61\\
 \textbf{8}	&22   &        18   &         8  &         13  &         69  &         66&           64    &        3   &         0  &          7   &        64\\   
 \textbf{9}	&24 & 	         20    &       12   &        16 &          71   &        69&           66  &          5  &          7 &           0  &         66\\   
\textbf{10}	& 60       &    55    &       53   &        52    &        6 &           3&            2   &        61     &      64 &          66 &           0\\  
\end{tabular}
\end{table}

\begin{figure}[ht!]
     \centering
     \begin{subfigure}[b]{0.45\textwidth}
         \centering
         \includegraphics[width=0.7\textwidth]{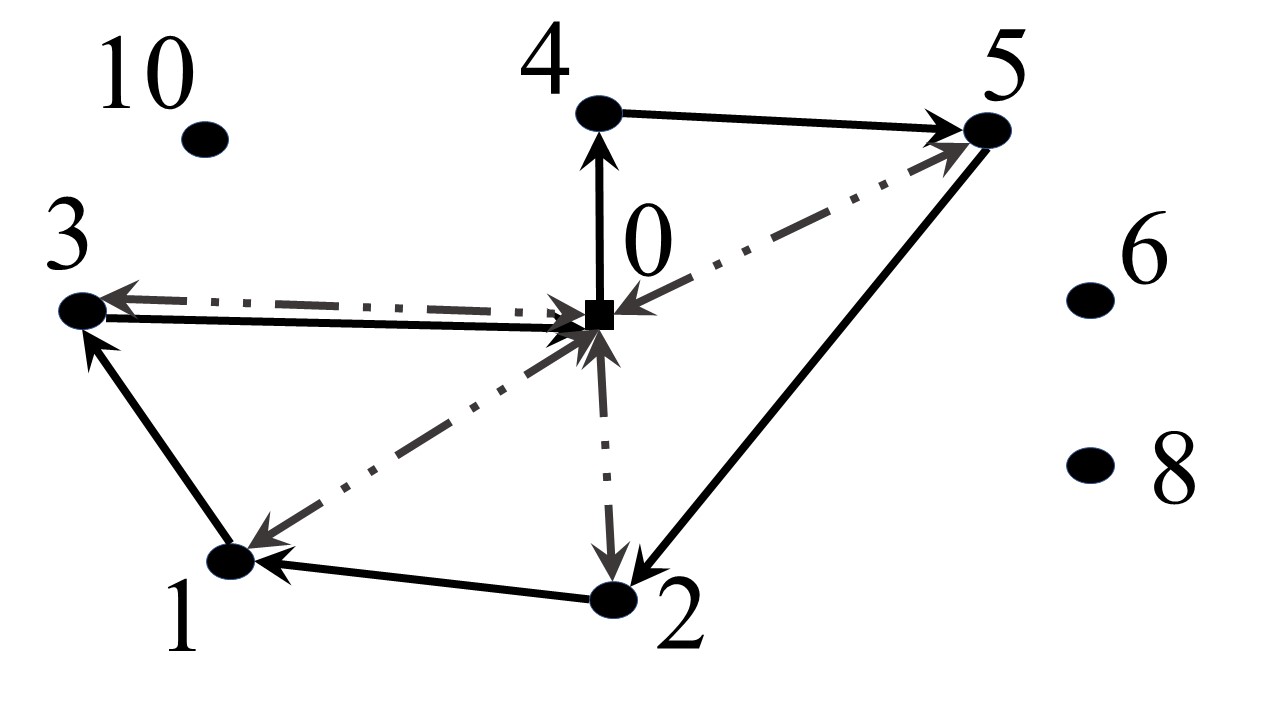}
         \caption{One truck.}
         \label{fig:onetruck}
     \end{subfigure}
     \hfill
     \begin{subfigure}[b]{0.45\textwidth}
         \centering
         \includegraphics[width=0.7\textwidth]{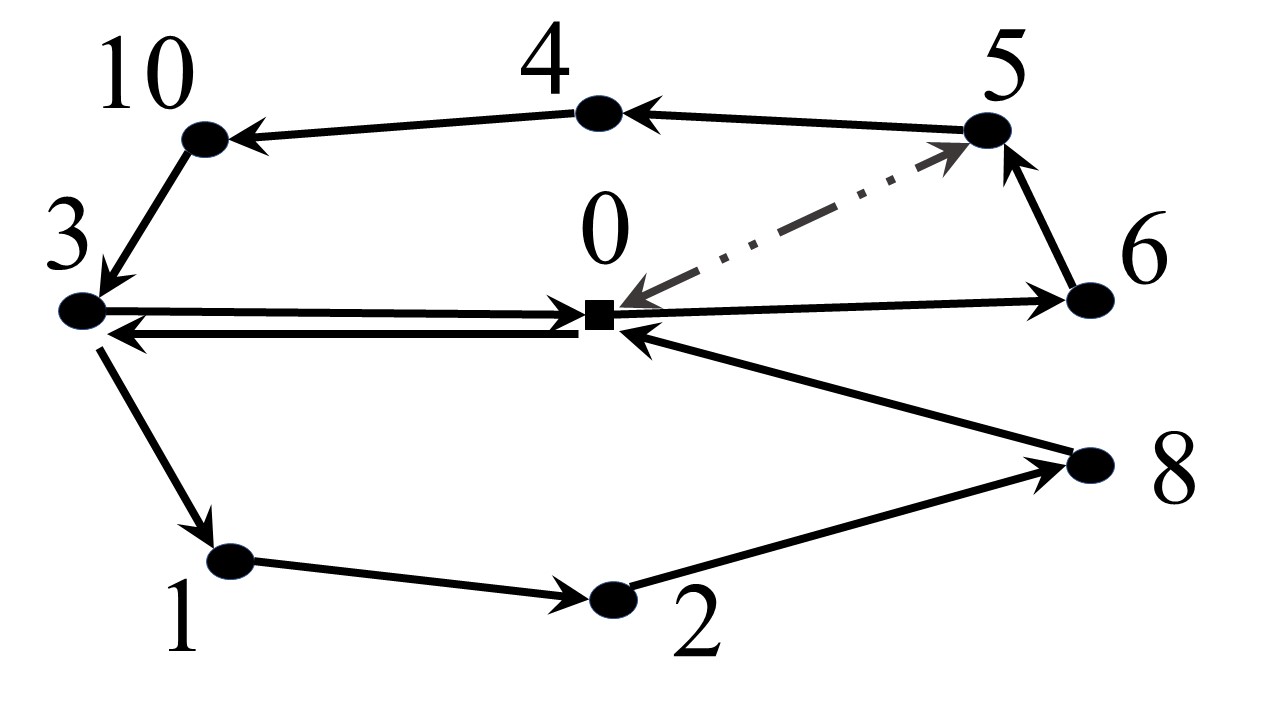}
         \caption{Two trucks.}
         \label{fig:twotruck}
     \end{subfigure}
     \hfill
       \caption{Optimal routes of the stochastic real-world example for two fleet sizes.}
        \label{fig:st}
\end{figure}

We begin by examining the stochastic case with $\omega = 0.8$. We analyze two situations separately: one in which we use a single truck and another in which we have two trucks that are assumed to operate only one route. Figure~\ref{fig:st} illustrates the optimal routes obtained in each of these situations. Table~\ref{tab:exampler} presents, for each situation, the load allocated to each truck compartment at the start of the route, along with the corresponding customer for each load. Note that, if $w=0.8$, some customers might not be visited. Finally, Table~\ref{tab:exampler2} (first two rows) displays the value of the objective function, the expected total distance to be traveled, the load transported (excluding round trips between a customer and the cooperative headquarters), the execution time of the solver, and the relative gap.

\begin{table}[ht!]
\centering
\caption{Compartment allocation in the stochastic real-world example.}\label{tab:exampler}
\begin{tabular}{c|cccc}
\toprule
              & \multicolumn{2}{c}{One truck} & \multicolumn{2}{c}{Two trucks} \\
Compartment  &   \multicolumn{2}{c}{\hrulefill} & \multicolumn{2}{c}{\hrulefill}   \\
              &  Customer    & Load               & Customer     & Load                  \\
\midrule 
   1 (truck 1)          &     4       &       3100            &      4      &         3100            \\
      2 (truck 1)       &     1       &       3000            &      3      &         3000             \\
         3 (truck 1)          &     3       &       1700            &      6     &         1700          \\
            4 (truck 1)          &     5       &       4500            &      5      &         4500           \\
 5 (truck 1)          &     2       &       3000            &      10      &         3000           \\
 \midrule
  1 (truck 2)          &     --       &       --            &      1      &         3500            \\
      2 (truck 2)          &     --       &       --            &      3      &         2959             \\
         3 (truck 2)          &     --       &       --            &      2     &         1541          \\
            4 (truck 2)          &     --       &       --            &      2      &         4500           \\
 5 (truck 2)          &     --       &       --            &      8      &         2799           \\   
 \botrule
\end{tabular}
\end{table}

\begin{table}[ht!]
\centering
\caption{Additional results in the real-world example.}
\label{tab:exampler2}
\begin{tabular}{l|ccccc}
\toprule
Scenario &  Objective    & Expected           & Total   & Time & Relative  \\
                   &               &  overall distance  &   load          &  (seconds) & gap \\
\midrule 
Stochastic,  $|K|=1$ & -2824 &   295  & 15,300 & 2.66 & 0 \\[0.2cm]
Stochastic, $|K|=2$ & -5927 &   241  & 30,599 & 300 & 0.008 \\[0.2cm]
Deterministic, $|K|=1$   & -2779.2 &   351  & 15,300 & 2.84 & 0 \\[0.2cm]
Deterministic, $|K|=2$ & -5959.2 &   201  & 30,600 & 102.16 & 0 \\
\botrule
 \end{tabular}
\end{table}

\begin{figure}[ht!]
     \centering
     \begin{subfigure}[b]{0.45\textwidth}
         \centering
         \includegraphics[width=0.7\textwidth]{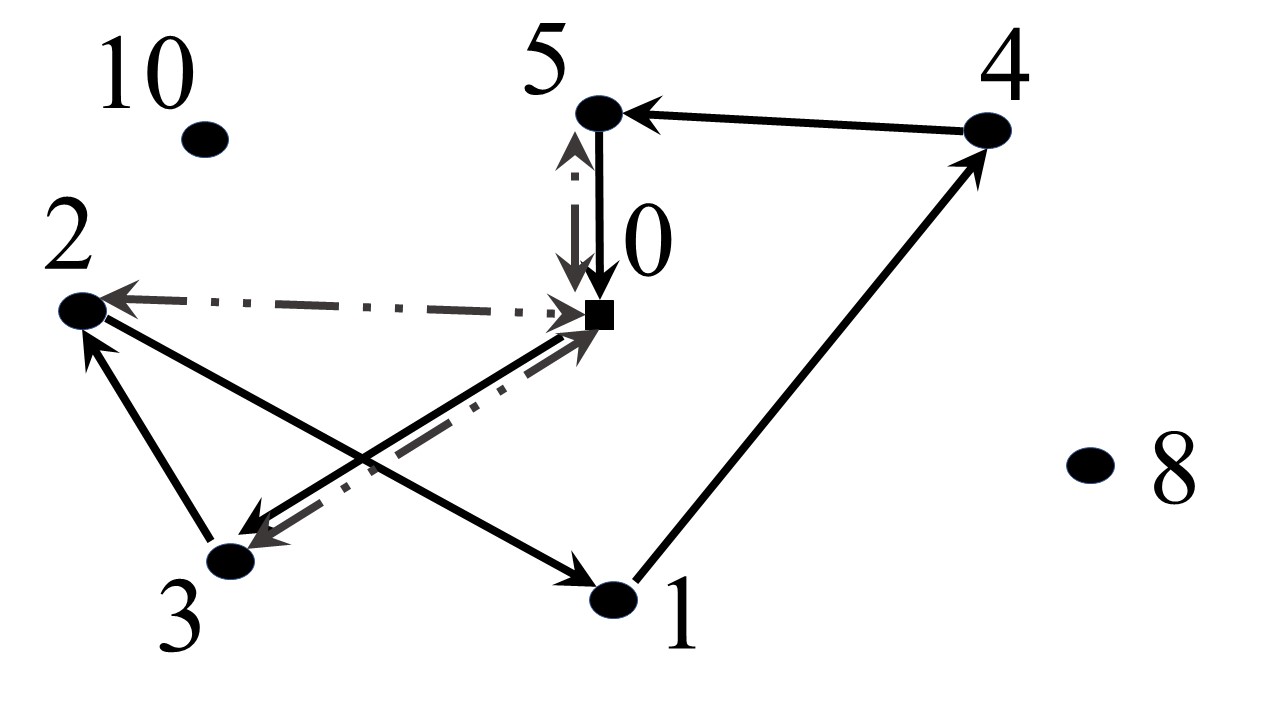}
         \caption{One truck.}
         \label{fig:onetruckd}
     \end{subfigure}
     \hfill
     \begin{subfigure}[b]{0.45\textwidth}
         \centering
         \includegraphics[width=0.7\textwidth]{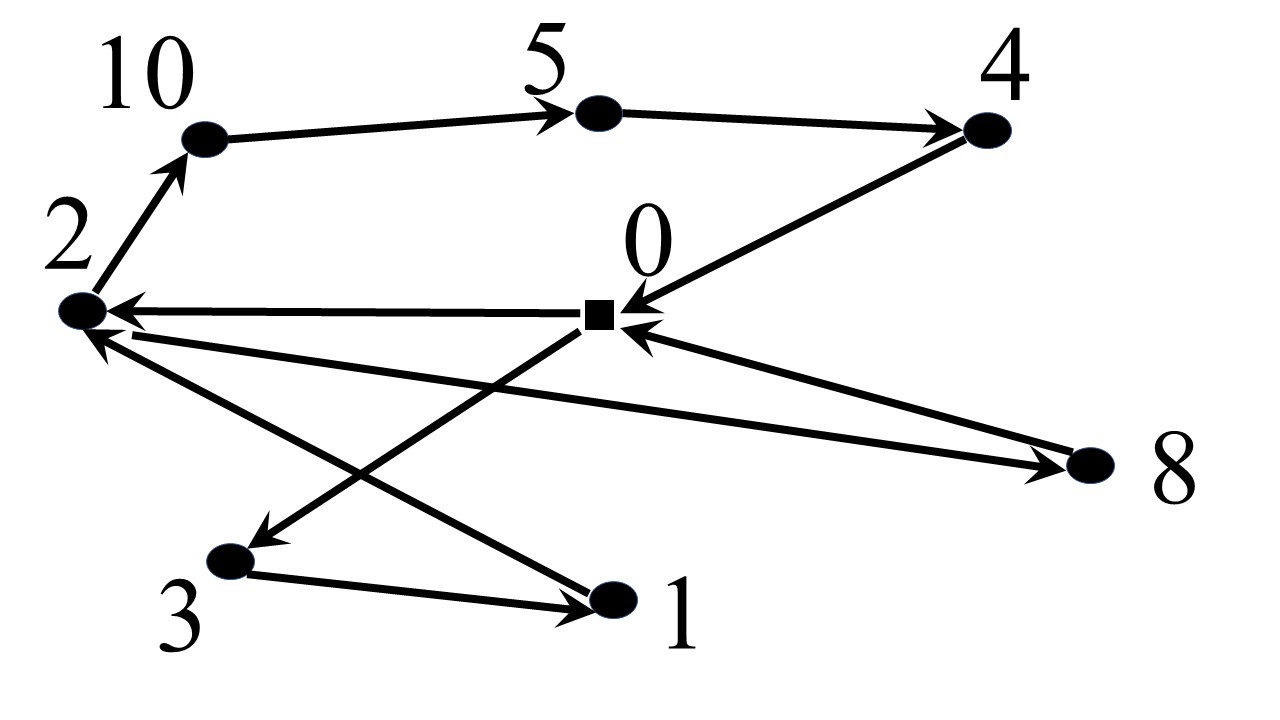}
         \caption{Two trucks.}
         \label{fig:twotruckd}
     \end{subfigure}
     \hfill
       \caption{Optimal routes of the deterministic real-world example for two fleet sizes.}
        \label{fig:dt}
\end{figure}

Next, we examine the deterministic case using the same weighting factor as before, $\omega=0.8$, and considering both one and two trucks. Figure~\ref{fig:dt} illustrates the optimal routes obtained in these situations. Tables~\ref{tab:exampled} and  \ref{tab:exampler2} (last two rows) provide results for the deterministic case that are analogous to those presented in Tables~\ref{tab:exampler} and  \ref{tab:exampler2} (first two rows)  for the stochastic case. 

\begin{table}[ht!]
\centering
\caption{Compartment allocation in the deterministic real-world example.} \label{tab:exampled}
\begin{tabular}{c|cccc}
\toprule
              & \multicolumn{2}{c}{One truck} & \multicolumn{2}{c}{Two trucks} \\
Compartment  &   \multicolumn{2}{c}{\hrulefill} & \multicolumn{2}{c}{\hrulefill}   \\
              &  Customer    & Load               & Customer     & Load                  \\
\midrule 
   1 (truck 1)          &     1       &       3300            &      5      &         4000            \\
      2 (truck 1)       &     4       &       2951            &      4      &         2951             \\
         3 (truck 1)          &     2       &       1549            &      5     &         885          \\
            4 (truck 1)          &     3       &       4500            &      2      &         4464           \\
 5 (truck 1)          &     5       &       3000            &      10      &         3000           \\
 \midrule
  1 (truck 2)          &     --       &       --            &      1      &         3300            \\
      2 (truck 2)          &     --       &       --            &      3      &         3000             \\
         3 (truck 2)          &     --       &       --            &      2     &         1577          \\
            4 (truck 2)          &     --       &       --            &      8      &         4464           \\
 5 (truck 2)          &     --       &       --            &      3      &         2959           \\   
 \botrule
\end{tabular}
\end{table}

By comparing the results of the deterministic and stochastic cases, we observe that using two trucks is advantageous in both scenarios, as it results in reduced distance traveled and increased load transported. The stochastic objective value is only $5.4\%$ worse (higher), primarily because there is a $1/3$ probability that a round trip between customer~5 and the cooperative headquarters will be required. If this round trip is not needed, the distances and total loads are nearly identical in both cases, allowing for the service of an additional customer (customer~6) in the stochastic scenario. We therefore believe that the stochastic model can be a valuable tool for the company's technicians in charge of managing the truck routes. It is also worth noting that uncertainty in demand increases computation time; for instance, if two trucks are used in the stochastic case, 300 seconds are insufficient to find a solution with a full guarantee of optimality.

\subsubsection{Extent of the exact solving of the model}\label{sec:extent_solving}

Table~\ref{tab2} provides the relative gap and computation time for the solving of 24 generated problems, designed to evaluate the model's capability in providing exact solutions. These problems consider different values for the number of customers, number of vehicles, number of urgent customers, and weights in the objective function. We set a time limit of 10 hours. It can be observed that problems involving 9 or 10 customers present a positive gap, with computation times reaching up to 10 hours. Additionally, both the gap and computation time increase, especially, the higher the weight assigned  in the objective function to the transported load. This behavior justifies the use of heuristic algorithms to address real-world problems of medium to large scale.

\begin{table}[ht!]
\centering
\caption{Relative gaps and computation times for the exact solution of 24 test cases.\label{tab2}}
\begin{tabular}{ccccccc}
\toprule
Case&Number of	& Urgent	& Number of& $\omega$& Relative& Time\\
number& customers	& customers	& vehicles & & gap&  (seconds)\\
\midrule 
1 &4  &  2  & 1 & 0.2  &  0    & 0.062\\
2&4  &  2  & 1 & 0.8  &  0     & 0.021\\
3&5  &  3  & 1 & 0.2  &  0     & 0.391\\
4&5  &  3  & 1 & 0.8  &  0     & 0.141\\
5&6  &  3  & 2 & 0.2  &  0     & 42.375\\
6&6  &  3  & 2 & 0.8  &  0     & 0.146\\
7&7  &  4  & 2 & 0.2  &  0.131 & 60\\
8&7  &  4  & 2 & 0.2  &  0     & 89.125\\
9&7  &  4  & 2 & 0.8  &  0     & 0.36\\
10&8  &  4  & 2 & 0.2  &  0.232 & 60\\
11&8  &  4  & 2 & 0.2  &  0     & 153.094\\
12&8  &  4  & 2 & 0.8  &  0     & 0.359\\
13&9  &  5  & 2 & 0.2  &  0.318 & 60\\
14&9  &  5  & 2 & 0.2  &  0.299 & 300\\
15&9  &  5  & 2 & 0.2  &  0.289 & 600\\
16&9  &  5  & 2 & 0.2  &  0.269 & 3600\\
17&9  &  5  & 2 & 0.8  &  0     & 3.5\\
18&10 &  5  & 2 & 0.2  &  0.393 & 60\\
19&10 &  5  & 2 & 0.2  &  0.378 & 300\\
20&10 &  5  & 2 & 0.2  &  0.373 & 600\\
21&10 &  5  & 2 & 0.2  &  0.339 & 3600\\
22&10 &  5  & 2 & 0.2  &  0.311 & 18000\\
23&10 &  5  & 2 & 0.2  &  0.301 & 36000\\
24&10 &  5  & 2 & 0.8  &  0     & 3.797\\
\botrule
\end{tabular}
\end{table}

\subsection{Results obtained with the algorithms}\label{sec:heuristic_results}

\subsubsection{Solving of the small real-world example}\label{sec:heuristic_small} 

\begin{figure}[t!]
     \centering
     \begin{subfigure}[b]{0.45\textwidth}
         \centering
         \includegraphics[width=0.7\textwidth]{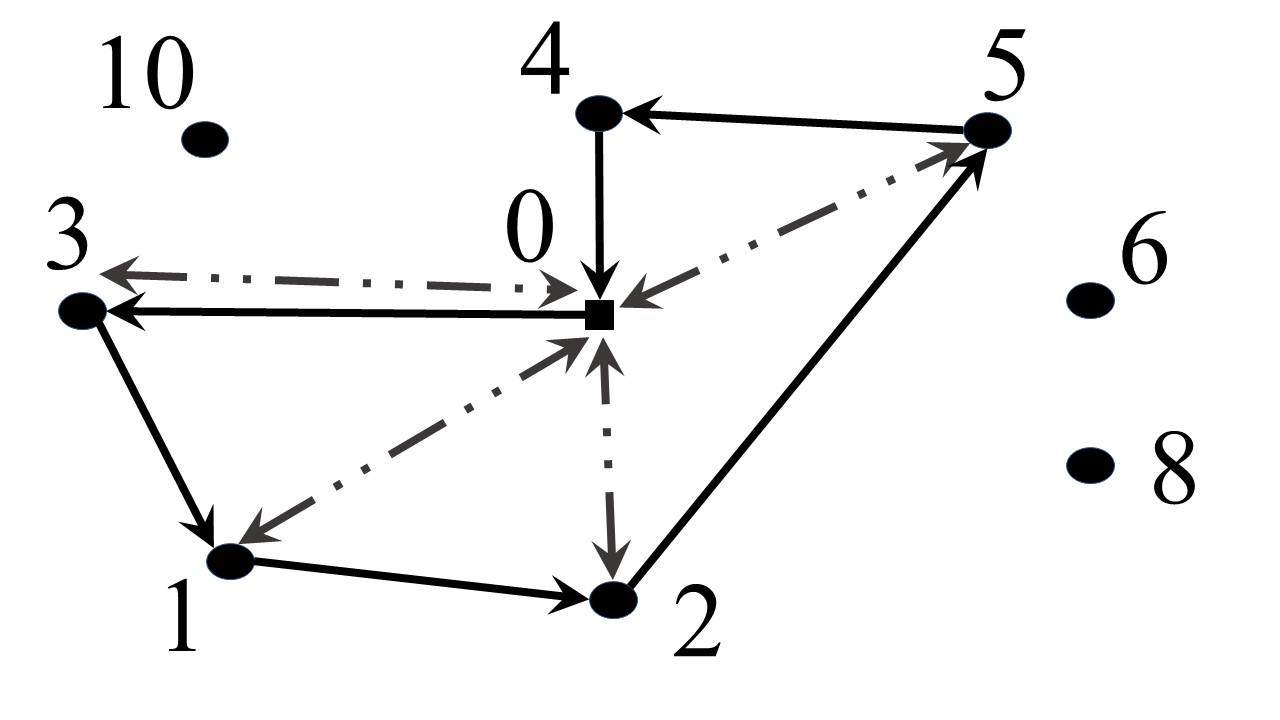}
         \caption{Stochastic case, one truck.}
         \label{fig:alreal1}
     \end{subfigure}
     \hfill
     \begin{subfigure}[b]{0.45\textwidth}
         \centering
         \includegraphics[width=0.7\textwidth]{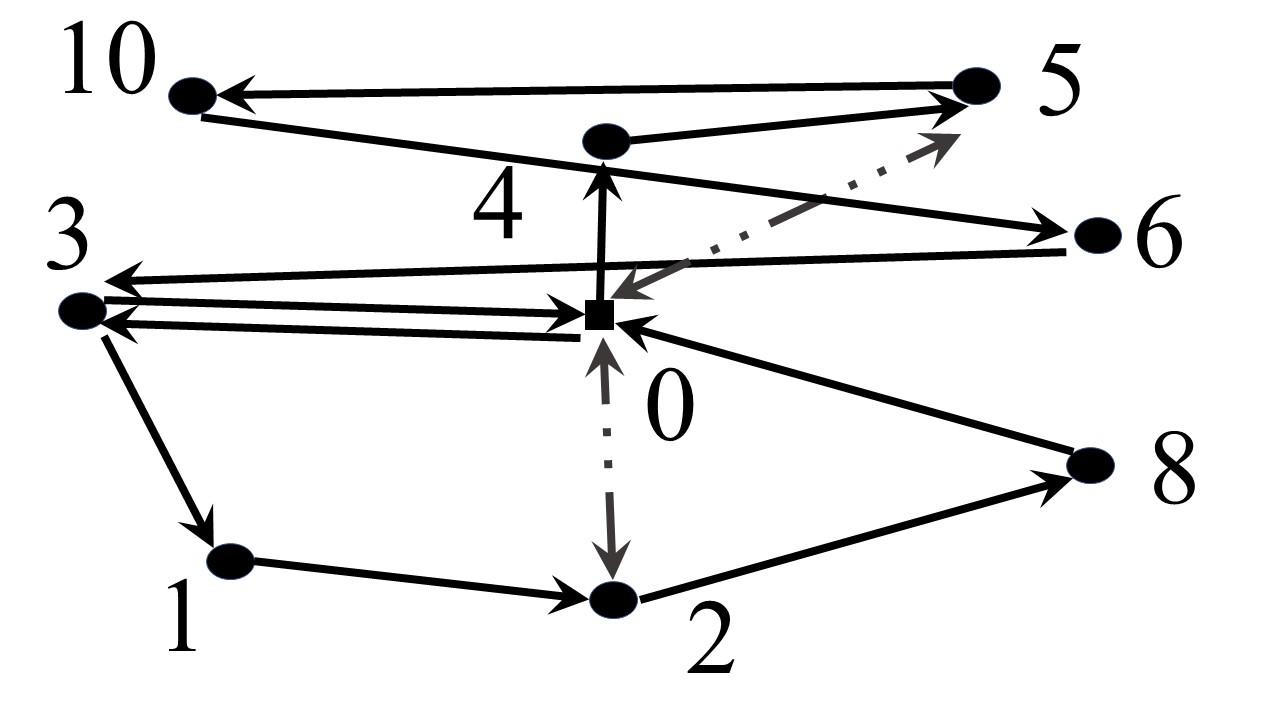}
         \caption{Stochastic case, two trucks.}
         \label{fig:alreal2}
     \end{subfigure}
     \hfill
     \begin{subfigure}[b]{0.45\textwidth}
         \centering
         \includegraphics[width=0.7\textwidth]{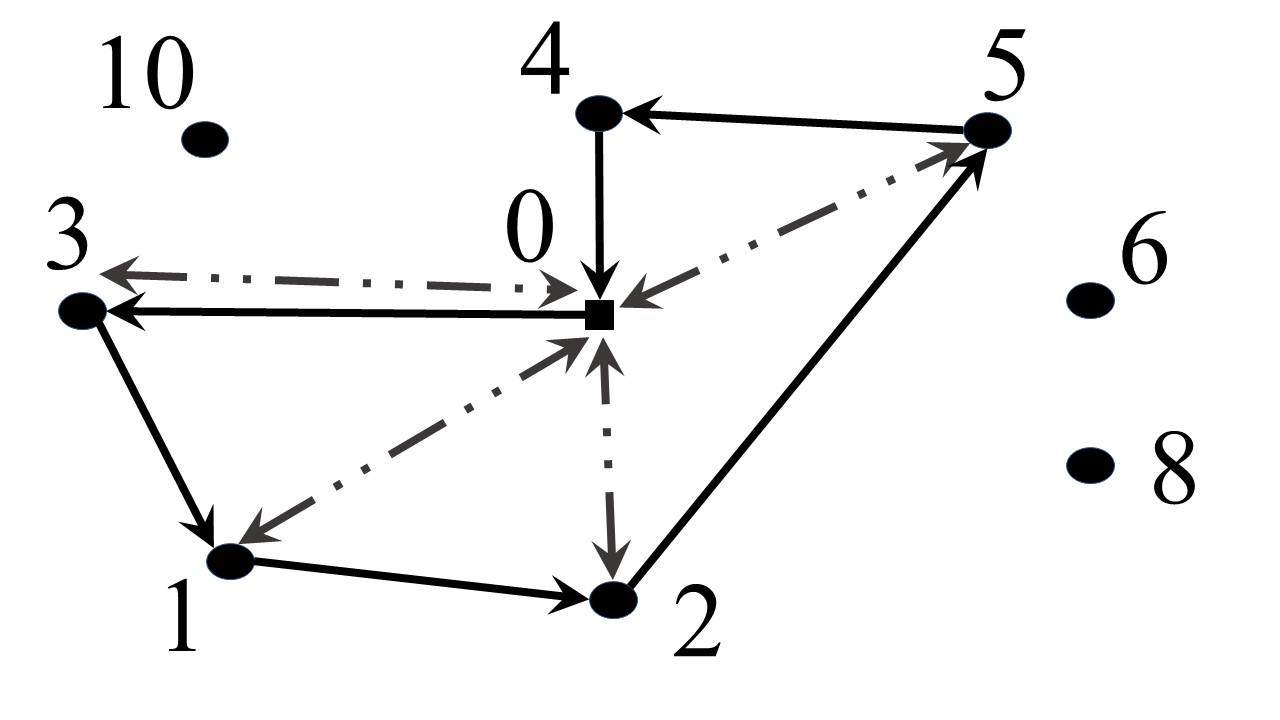}
         \caption{Deterministic case, one truck.}
         \label{fig:alreal3}
     \end{subfigure}
     \hfill
     \begin{subfigure}[b]{0.45\textwidth}
         \centering
         \includegraphics[width=0.7\textwidth]{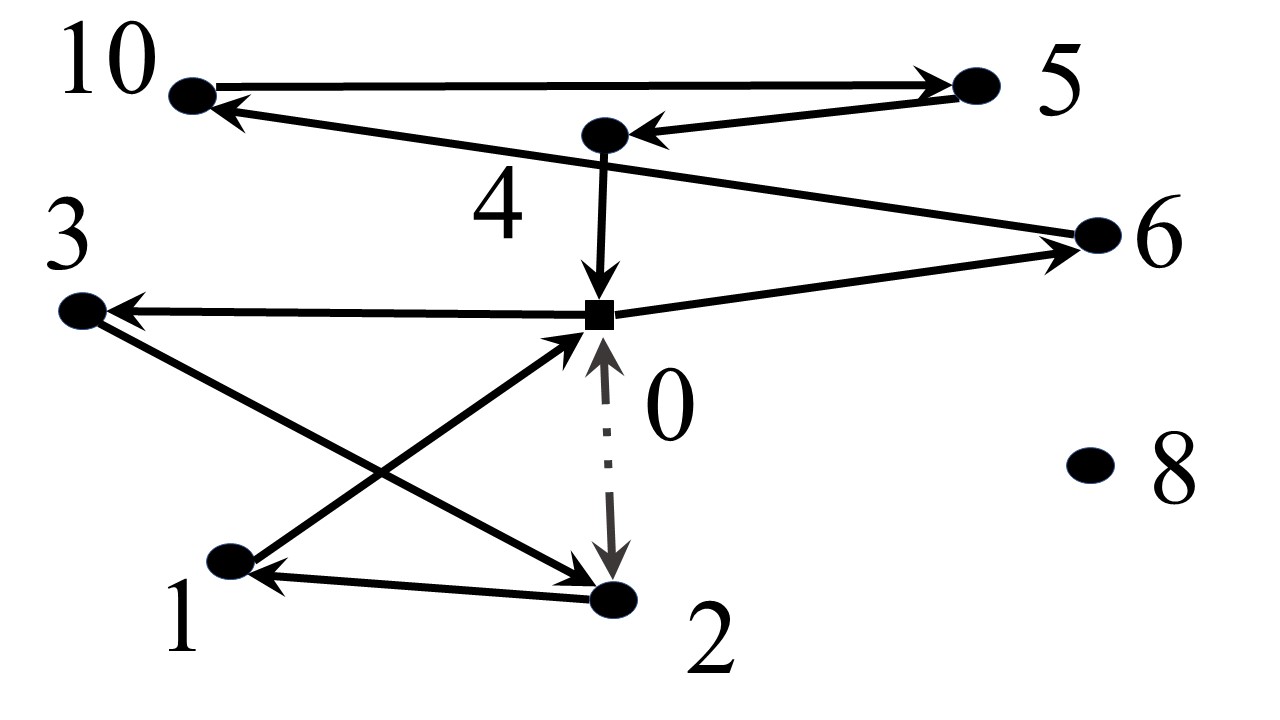}
         \caption{Deterministic case, two trucks.}
         \label{fig:alreal4}
     \end{subfigure}
     \hfill
       \caption{Route configurations in the real-world example obtained with the two-stage heuristic algorithm for different scenarios.}
        \label{fig:alreal}
\end{figure}

We have re-evaluated the small real-world example solved exactly and applied our algorithm to it. Figure~\ref{fig:alreal} and Table~\ref{tab:rexamplealg} present the results. Comparing the outcomes from the heuristic algorithm to those obtained from the exact solution of the model, we observe that the traveled distances match in some instances. However, in other cases, the distance can increase by as much as $8.95\%$. A similar pattern is noted for the transported load, which can decrease by up to $12.43\%$ in the worst-case scenario. Regarding the objective value, it either remains unchanged or, in the worst case, deteriorates by up to $13\%$. Nonetheless, the computation times for the heuristic algorithm are reduced by as much as $97\%$.

\begin{table}[ht!]
\centering
\caption{Results obtained with the two-stage heuristic algorithm in the real-world example.}
\label{tab:rexamplealg}
\begin{tabular}{l|cccc}
\toprule
Scenario           &  Objective    & Expected           & Total load  & Time      \\
               &               &  overall distance  &             &  (seconds) \\
\midrule 
Stochastic,    $|K|=1$ &    -2824           &     295               &     15,300        & 0.06  \\[0.2cm]
Stochastic,   $|K|=2$ &     -5352.46   &  254.33            & 27,779.5     & 19.18   \\[0.2cm]
Deterministic, $|K|=1$ & -2509.4        &  351               &  13,951      &  0.16 \\[0.2cm]
Deterministic,  $|K|=2$ &   -5183.8            &    219                &    26,795 & 18.06  \\
\botrule
\end{tabular}
\end{table}

\subsubsection{Set of designed instances}\label{sec:new_instances}

Since the MC-VRPSD considered in this work is a novel problem in the literature, no  established baseline exits for evaluating the proposed solving algorithms. To fill this gap, we generated three new problem sets with varying characteristics, resulting in a total of 35 instances. These instances are adaptations of 14 problems proposed by \cite{elf08} for the MC-VRP. The original instances, created by \cite{Chr79} for the VRP, include between $50$ and $199$ customers. For further details about the nature of these problems, we recommend referring to the original work of \cite{elf08}.

We classified these instances into three sets according to the statistical distribution of demand, the number of trucks, the number and capacity of truck compartments, and the number of urgent customers or customers with accessibility restrictions for one of the trucks:

\begin{itemize}
    \item Set 1: All customers have urgent demands. 
    \item Set 2: The first half of the customers have urgent demands, while the rest are non-urgent. 
    \item Set 3: The first half of the customers are urgent, and the remainder are non-urgent. Two vehicles are used, but the vehicle with 6 compartments is restricted from accessing $20\%$ of the customers. 
\end{itemize}

In all three cases, the stochastic demand follows a normal distribution with a mean equal to the demand requested in the deterministic case (as given in the original instances) and a standard deviation set at  $30\%$ of that demand. The total vehicle capacity is the sum of its compartment capacities. Additional relevant data are provided in Table~\ref{tab:instances}, showing variety in the number, urgency, and accessibility constraints of customers, as well as in the compartment configuration.


\begin{table}[ht!]
\centering
\caption{Characteristics of the generated 35 MC-VRPSD test problems.}
\label{tab:instances}
\begin{tabular*}{\textwidth}{@{\extracolsep\fill}l|ccccc}
\toprule%
          & \multicolumn{3}{@{}c@{}}{Customers}  & \multicolumn{2}{@{}c@{}}{Trucks} \\
Problem   & \multicolumn{3}{c}{\hrulefill} & \multicolumn{2}{c}{\hrulefill}  \\
 number   & Number  & Urgent & Accessibility  & Number of & Compartments     \\
          &         &        & restrictions  &  compart.        & capacity    \\
           \midrule
Set 1     &         &        &             &          &               \\
vrnpc1    &    50   &  50    &     --      &    5      &   45, 35, 30, 30, 20      \\ 
vrnpc1b    &    50   &  50    &     --      &    6      &   40, 35, 30, 25, 20, 10      \\
vrnpc2    &    75  &  75    &     --      &    5      &   40, 30, 30, 20, 20     \\ 
vrnpc2b    &    75   &  75    &     --      &    6      &   35, 30, 25, 20, 15, 15      \\
vrnpc3    &    100  & 100    &     --      &    5      &  50, 45, 40, 35, 30     \\ 
vrnpc3b    &   100    &  100     &     --      &    6      &   45, 40, 35, 30, 30, 20  \\
vrnpc4    &    150  & 150    &     --      &    5      &  50, 45, 40, 35, 30     \\
vrnpc4b    &    150  & 150    &     --      &    6      &  45, 40, 35, 30, 30, 20     \\
vrnpc5    &    199  & 199    &     --      &    5      &  50, 45, 40, 35, 30     \\
vrnpc5b    &    199  & 199    &     --      &    6      & 45, 40, 35, 30, 30, 20     \\
vrnpc11    &    120  & 120    &     --      &    5      &  50, 45, 40, 35, 30     \\
vrnpc11b    &    120  & 120    &     --      &    6      &  45, 40, 35, 30, 30, 20     \\
vrnpc12    &    100  & 100    &     --      &    5      &  50, 45, 40, 35, 30     \\
vrnpc12b    &    100  & 100    &     --      &    6      &  45, 40, 35, 30, 30, 20    \\[0.2cm]
Set 2     &         &        &             &          &                \\
vrnpc1    &    50   &  First 50$\%$    &     --      &    5      &   45, 35, 30, 30, 20      \\ 
vrnpc1b    &    50   &  First 50$\%$    &     --      &    6      &   40, 35, 30, 25, 20, 10      \\
vrnpc2    &    75  &  First 50$\%$    &     --      &    5      &   40, 30, 30, 20, 20     \\ 
vrnpc2b    &    75   &  First 50$\%$    &     --      &    6      &   35, 30, 25, 20, 15, 15      \\
vrnpc3    &    100  & First 50$\%$    &     --      &    5      &  50, 45, 40, 35, 30     \\ 
vrnpc3b    &   100    &  First 50$\%$     &     --      &    6      &   45, 40, 35, 30, 30, 20  \\
vrnpc4    &    150  & First 50$\%$     &     --      &    5      &  50, 45, 40, 35, 30     \\
vrnpc4b    &    150  & First 50$\%$   &     --      &    6      &  45, 40, 35, 30, 30, 20     \\
vrnpc5    &    199  & First 50$\%$    &     --      &    5      &  50, 45, 40, 35, 30     \\
vrnpc5b    &    199  & First 50$\%$    &     --      &    6      & 45, 40, 35, 30, 30, 20     \\
vrnpc11    &    120  & First 50$\%$    &     --      &    5      &  50, 45, 40, 35, 30     \\
vrnpc11b    &    120  & First 50$\%$    &     --      &    6      &  45, 40, 35, 30, 30, 20     \\
vrnpc12    &    100  & First 50$\%$    &     --      &    5      &  50, 45, 40, 35, 30     \\
vrnpc12b    &    100  & First 50$\%$   &     --      &    6      &  45, 40, 35, 30, 30, 20    \\[0.2cm]
Set 3     &         &        &             &          &               \\
vrnpc1    &    50   &  First 50$\%$    &     4, 39, 1, 34, 23,      &    5, 6      &   45, 35, 30, 30, 20      \\ 
    &       &      &    43, 14, 18, 33, 21     &          &   40, 35, 30, 25, 20, 10      \\
vrnpc2    &    75  &  First 50$\%$    &     68, 39, 1, 34, 43, 14, 59, 51,      &    5, 6&  40, 30, 30, 20, 20     \\     &      &      &    21, 54, 7, 9, 15, 67, 37      &          &   35, 30, 25, 20, 15, 15      \\
vrnpc3    &    100  & First 50$\%$    &68, 39, 1, 34, 87, 43, 14, 82, 59, 51, &    5, 6    &  50, 45, 40, 35, 30     \\ 
    &       &       &  85, 21, 54, 74, 7, 73, 79, 37, 83, 97 &          &   45, 40, 35, 30, 30, 20  \\
vrnpc4    &    150  & First 50$\%$     &68, 129, 43, 14, 51, 85, 21, 106, 74, 7, &  5, 6    &  50, 45, 40, 35, 30     \\
    &      &    &     73, 79, 37, 105, 110, 34, 143, 126, 89, 33,      &          &  45, 40, 35, 30, 30, 20     \\
        &      &    &    84, 70, 142, 42, 38, 11, 20, 28, 124, 44      &          &       \\
vrnpc5    &    199  & First 50$\%$    & 68, 87, 1, 67, 129, 162, 43, 14, 187, 51, & 5, 6   &  50, 45, 40, 35, 30     \\
    &  &  &  85, 21, 106, 182, 74, 7, 73, 79, 37, 105,   &          & 45, 40, 35, 30, 30, 20     \\
       &  &  &  110, 165, 34, 189, 126, 89, 172, 33, 84, 163,   &          &       \\
          &  &  &  70, 193, 42, 166, 11, 148, 156, 20, 44, 121   &          &       \\
vrnpc11    &    120  & First 50$\%$    &  68, 39, 1, 34, 87, 43, 14, 82,  &    5, 6      &  50, 45, 40, 35, 30     \\
  &      &      &    59, 51, 97, 85, 21, 106, 54, 74, &          &  45, 40, 35, 30, 30, 20     \\
    &      &      &    7, 73, 79, 109, 37, 89, 100, 117 &          &      \\
vrnpc12    &100  & First 50$\%$    & 68, 39, 1, 34, 87, 43, 14, 82, 59, 51, &    5, 6      &  50, 45, 40, 35, 30     \\
&  &   &  85, 21, 54, 74, 7, 73, 79, 37, 83, 97  &          &  45, 40, 35, 30, 30, 20    \\
\botrule
\end{tabular*}
\end{table}


It should be noted that when initializing the algorithms, we set the value of $m=5000$ for the maximum number of iterations in the tabu search algorithm, and $\kappa=1$ for the maximum number of possible exchanges between two routes, as the number of customers visited on each route is relatively small.

\subsubsection{Results on sets of instances} \label{sec:heuristic_results_instances}

\noindent \textbf{Experiment 1: Performance on MC-VRP instances and comparison with literature results.} To obtain a first validation of the quality of our algorithms, we analyzed their performance on MC-VRP instances. This allows for a comparison between the solutions generated by our method and those reported by \cite{elf08} and \cite{Men10}.

Table~\ref{tabcomp1} summarizes the results obtained for both phases of the heuristic. The first column indicates the problem number, while columns~2 and 3 display the objective values reported by \cite{Men10} and \cite{elf08}, respectively. Columns~ 4 and 5 present the results provided by our constructive algorithm and our ITS. It is important to note that we considered the objectives reported by \cite{Men10} when a memetic algorithm is implemented and those from \cite{elf08} when a tabu search is implemented. Regarding the solutions produced by the constructive stage, our savings-based heuristic successfully generated feasible routes. The ITS improves upon the objectives reported by \cite{Men10} or \cite{elf08} in five instances, closely matches them in most others, and outperforms both in two problems. While our heuristic is efficient for solving the MC-VRP, it should be noted that it was  specifically designed for the MC-VRPSD. This comparison with existing MC-VRP solutions serves to illustrate that our algorithm performs well on problems that may be similar.\\

\begin{table}[ht!]
\centering
\caption{Computational results for 14 test instances of the MC-VRP.}
\label{tabcomp1}
\begin{tabular}{l|cccc}
\toprule
Problem number & \cite{Men10} & \cite{elf08} & Our CW & Our ITS   \\
\midrule
 vrpnc1        & 	{\bf 524.6}  & {\bf 524.6}	&625.55 &  	532.02\\
 vrpnc2        & 	857.3  & {\bf 850.0}	&997.74& 	855.98\\
 vrpnc3        & 	841.6  & {\bf 831.3}	&1001.05& 	868.40\\
 vrpnc4        & 	{\bf 1045.9} & 1061.1	&1276.52& 	1073.06\\
 vrpnc5        & 	1381.9 & {\bf 1348.3}	&1643.49& 	1361.12\\
 vrpnc6        & 	{\bf 556.7}  & 575.9	&652.01& 	560.88\\
 vrpnc7	       &    {\bf 924.3}  &  970.8  &1058.42& 	936.36\\
 vrpnc8        & 	{\bf 877.6}  & 888.6	&1053.74	& 894.30\\
 vrpnc9        & 	1223.1 & 1232.1	&1399.72& 	{\bf 1218.64}\\
 vrpnc10       & 	{\bf 1497.5} & 1538.6	&1678.67& 	1522.51\\
 vrpnc11       & 	1046.9 & {\bf 1043.3}	&1200.95& 	1060.43\\
 vrpnc12  	   &    821.0  & {\bf 819.5}	&939.99& 	843.25\\
 vrpnc13       & 	1586.9 & 1582.2	&2368.85& 	{\bf 1522.61}\\
 vrpnc14       &    {\bf 868.1}  & 868.6	&945.57& 	940.08\\
 \botrule
\end{tabular}
\end{table}

\noindent \textbf{Experiment 2: Performance on new MC-VRPSD instances.}  
Focusing on the specific problem of interest, Table~\ref{tabcomp2} presents the performance of the proposed heuristic on the 35 MC-VRPSD instances described in Table~\ref{tab:instances}, tailored to the characteristics of our problem. The first column indicates the problem number, followed by three columns detailing the fixed distance, the total expected distance, and the computation time (in seconds) for the constructive phase. The subsequent columns refer to the ITS solution, with columns~5 to 8 displaying the fixed distance and the expected total distance of the improvement phase (reporting the best value over five runs), the total computation time (in seconds) including both the constructive and improvement stages, and the number of routes created. Note that it is assumed that different products and products for different customers cannot be mixed, which leads to an increased number of routes compared to scenarios where these assumptions are not made, as well as a lower occupancy or filling of the compartments. The last column shows the occupancy rate of the vehicles, calculated as the overall average of the average occupancy per route. These occupancy rates range from $33\%$ to $87\%$. The lowest rates are found in Set 3, always below $60\%$, while the highest rates are achieved in Set 1, always above $50\%$ except for problem vrpnc11. Furthermore, the lowest rates are consistently found with 120 customers, whereas the highest rates occur with 75 customers.

\begin{table}[ht!]
\centering
\caption{Computational results for the generated 35 MC-VRPSD test problemns.}
\label{tabcomp2}
\begin{tabular}{l|cccccccc}
\toprule
 & \multicolumn{3}{c}{CW} & \multicolumn{5}{c}{ITS} \\
Problem & \multicolumn{3}{c}{\hrulefill} & \multicolumn{5}{c}{\hrulefill} \\
number  &   Fixed    & Expected   & Time & Fixed    & Expected  & Total & Number of  & Occupancy   \\
   &   distance & distance   &   (seconds)   & distance & distance  & time (s) & routes & rate ($\%$)  \\
\midrule
Set 1    &            &            &      &          &           &     &  &   \\
vrpnc1   &   856.99   &  981.92    & 0.05 &  751.65  & 874.00    &  110   & 11 &67.39 \\
vrpnc1b  &   833.83   &  994.37    & 0.05  &   742.11 &    867.05  &  456   &10&76.48 \\ 
vrpnc2    &   1352.80     & 1568.57   &0.11 &   1216.97  &  1428.44   &     260   &18&80.60  \\
vrpnc2b   &  1277.90   &   1613.63  & 0.14 &1179.15  &   1434.62        &    483   &17&86.81   \\
vrpnc3    &   1594.63 & 1844.10  & 0.16 &  1449.37   & 1698.84          &    159  &21&52.42    \\
vrpnc3b    & 1436.62  & 1694.48 & 0.18 & 1308.22    & 1557.79          &  178   &18&61.35     \\
vrpnc4    & 2119.34  & 2491.19   &0.40 & 1989.58    &  2357.61         &     177  &31&54.12   \\
vrpnc4b  & 1921.17    &  2309.80  & 0.45     & 1771.66    &   2140.08 &       340 &26&63.59  \\
vrpnc5    & 2709.64  &3190.06   & 0.92   & 2516.97   & 2997.36  & 218   &41&58.41      \\
vrpnc5b  &2417.45   & 2960.43  &1.03 &2249.09    & 2729.75     & 420    &35&67.47     \\
vrpnc11  & 3099.82  & 3711.78  &0.25    &2758.42   &  3370.38         & 168     &25&43.81    \\
vrpnc11b &  2558.73  &  3201.77  &  0.28 &  2442.20   &  3024.16 &   181  &21&54.11 \\
vrpnc12  & 1605.35   &  1893.30  & 0.15  &  1531.06   &  1819.61 &    250   & 22&62.50  \\
vrpnc12b &  1452.24  &  1787.99  &  0.17    &   1363.20   &  1663.72 &  305   & 19&74.20    \\[0.2cm]
Set 2    &            &            &      &          &           &     &     \\
vrpnc1  &   802.09  &  889.38  & 0.05  &  760.23 &   814.59    &  90    & 11&56.57   \\
vrpnc1b &   797.60  &  874.45   &  0.05    & 714.19  & 770.27          &  160  &10&66.01 \\
vrpnc2  & 1165.75 &  1302.85  &  0.10    & 1060.94         & 1147.06   &  199  &14&75.40  \\
vrpnc2b & 1167.85  & 1371.19 &  0.15    &  1092.17  &  1180.47     & 145     &  15&79.43  \\
vrpnc3 & 1536.92 & 1668.15 &  0.15  &  1444.20  & 1575.43          &  95     & 21&46.55  \\
vrpnc3b & 1395.49  & 1564.98  & 0.17     &  1284.45        &  1415.98         &  89    & 17&54.28\\
vrpnc4    & 2104.51  &  2300.67  &  0.39    & 1981.25  &   2177.47        &  175  &31&47.21      \\
vrpnc4b    & 1977.43  & 2173.59   & 0.41     &  1752.01   &  1948.47 &  377     & 26&54.94  \\
vrpnc5  & 2687.62 &  2944.24 & 0.89   & 2485.66  &  2735.13  &   206  & 41&49.08    \\
vrpnc5b & 2388.69  &  2661.97  & 0.95     & 2172.37   & 2422.13   &    290   & 34&57.95  \\
vrpnc11 &  3051.15 & 3483.26 &  0.22&  2773.60        &  3205.71         &   130    & 24&37.39  \\
vrpnc11b & 2619.38  & 3080.06 & 0.25     &  2480.30        &  2912.41  &  175  &  20&44.89    \\
vrpnc12  &  1519.60  & 1640.19   & 0.19     & 1461.18   & 1581.75   & 40   &  21&54.52     \\
vrpnc12b &   1407.81 &  1572.45   &  0.22    & 1355.78   &  1663.72         &  165  & 18&64.03 \\[0.2cm]
Set 3    &            &            &      &          &           &    &      \\
vrpnc1    & 890.52  &   962.48  &  0.10    &   734.53       &  792.70         & 98 &10&50.53  \\
vrpnc2    & 1226.37 &  1431.68   &  0.22    & 1112.36  & 1217.71  &  121   & 16&59.49    \\
vrpnc3    &  1547.10 &   1684.97 &  0.15    &  1397.85   &  1529.09         & 112  &19&41.91   \\
vrpnc4 & 2065.57  & 2290.76  &      0.78  & 1897.65          &2093.81      & 131 & 28&41.80 \\
vrpnc5    & 2644.59  &   2894.29 &  1.55    & 2444.27         & 2694.94    &    55  & 38&43.96   \\
vrpnc11  & 2971.21   &  3407.36          &  0.45    & 2740.78  &  3172.89    &   287  &23&33.55 \\
vrpnc12    &  1546.01  &  1677.53   & 0.32  &  1481.37  &   1601.94   &    121  &19 &48.88   \\
\botrule
\end{tabular}
\end{table}

\section{Conclusions}\label{sec:conclusions}

In this paper, we introduce a novel variant of vehicle routing problems characterized by a heterogeneous fleet comprising compartmentalized vehicles with varying capacities. These vehicles are responsible for distributing various products to a set of customers with stochastic demands. Notably, each compartment is restricted to containing only one type of product for a single customer. Additionally, some vehicles cannot access certain customers, and specific orders are classified as urgent. To optimize this routing problem, we propose a two-stage heuristic algorithm that begins with a constructive heuristic. In this initial phase, we generalize the savings-based algorithm of Clarke and Wright. The solution generated from this adaptation serves as the starting point for an iterated tabu search algorithm, specifically designed to handle stochastic demands and multi-compartment vehicles--an approach that, to the best of our knowledge, is novel in the existing literature. We illustrate the model and algorithms using both simulated data and a real-world case study. Furthermore, we compare the results of our algorithms with  well-known methods from the literature for particular problems and create a collection of benchmark instances for future research in this field.

The main contribution of this work is the proposal of a model that accurately captures real-world scenarios by incorporating several constraints in a novel manner. Specifically, the research conducted herein was motivated by the operations of an agricultural cooperative in a Spanish region of approximately 13,000 km$^2$, which hosts over 15 similar companies. The model we propose is versatile and applicable to other types of businesses, including those involved in milk collection, cereal distribution, and gasoline logistics. The algorithms developed can serve as valuable decision-making aids for optimizing transportation within these sectors.

These algorithms are implemented in the \texttt{R} programming language, enabling their integration into a comprehensive tool to assist company managers in decision-making processes. This tool can complement other systems, such as those used for consumption forecasting on farms or generating relevant reports. Furthermore, it is possible to develop user-friendly interfaces that simplify managerial tasks and facilitate simulations for advanced users using specific \texttt{R} libraries. Such tools represent a significant advantage for companies, as they automate certain tasks more effectively than traditional human-based methodologies, which tend to be less efficient, or existing commercial tools, which often lack transparency, are expensive, and difficult to adapt.

Regarding possible future research directions, it would be worthwhile to explore the adaptation of other established metaheuristics, such as simulated annealing or genetic algorithms, allowing for comparative analyses with our algorithms using the benchmarks developed in this paper. Additionally, investigating alternative approaches to address accessibility restrictions--such as utilizing both trucks and trailers when available (see \citealp{Der13} or \citealp{Dav23})--in stochastic scenarios could provide further insights and lead to reduced transportation costs.

\bmhead{Acknowledgements}
This work is part of the R+D+I projects PID2021-12403030NB-C31 and PID2021-124030NB-C32, granted by MICIU/AEI/10.13039/501100011033/ and by ``ERDF A way of making Europe''/EU. This research was also funded by {\it Grupos de Referencia Competitiva} ED431C 2021/24 from the \emph{Conseller{\'i}a de Cultura, Educaci{\'o}n e Universidade, Xunta de Galicia}.

\section*{Statements and Declarations}
\noindent \textbf{Competing Interests} The authors have no conflicts of interest to declare that are relevant to the content of this article.

\bibliography{biblio}

\end{document}